\documentclass{article}

\usepackage{graphicx}
\usepackage{amsmath}

\newtheorem{theorem}{Theorem}

\newtheorem{corollary}[theorem]{Corollary}

\newtheorem{definition}[theorem]{Definition}

\newtheorem{lemma}[theorem]{Lemma}
\newtheorem{notation}[theorem]{Notation}

\newtheorem{proposition}[theorem]{Proposition}
\newtheorem{remark}[theorem]{Remark}

\newenvironment{proof}[1][Proof]{\textbf{#1.} }{\ \rule{0.5em}{0.5em}}

\begin{document}

\title{Isomorphic groupoid $C^{\ast }$-algebras associated with different Haar
systems}
\author{M\u{a}d\u{a}lina Roxana Buneci\thanks{%
This work was partly supported by the MEC-CNCSIS grant At127/2004 and by
the Postdoctoral Training Program HPRN-CT-2002-0277.}
\\
University Constantin Br\^{a}ncu\c{s}i, T\^{a}rgu-Jiu}
\date{}
\maketitle

\begin{abstract}
We shall consider a locally compact groupoid endowed with a Haar system $\nu 
$ and having proper orbit space. We shall associated to each appropriated
cross section $\sigma :G^{\left( 0\right) }\rightarrow G^{F}$ for $%
d_{F}:G^{F}\rightarrow G^{\left( 0\right) }$ (where $F$ is a Borel subset of 
$G^{\left( 0\right) }$ meeting each orbit exactly once) a $C^{\ast }$%
-algebra $M_{\sigma }^{\ast }\left( G,\nu \right) $. We shall prove that the
\ $C^{\ast }$-algebras $M_{\sigma }^{\ast }\left( G,\nu _{i}\right) $ $i=1,2$
associated with different Haar systems are $\ast $-isomorphic
\end{abstract}

{\small AMS 2000 Subject Classification: 22A22, 43A22, 43A65, 46L99.}

{\small Key Words: locally compact groupoid, }$C^{\ast }$-algebra, $\ast $%
-isomorphism.

\section{Introduction}

The $C^{\ast }$-algebra of a locally compact groupoid was introduced by J.
Renault in \cite{re1}. The construction extends the case of a group: the
space of continuous functions with compact support on groupoid is made into
a $\ast $-algebra and endowed with the smallest $C^{\ast }$-norm making its
representations continuous. In order to define the convolution on groupoid
one needs to assume the existence of a Haar system which is an analogue of
Haar measure on a group. Unlike the case for groups, Haar systems need not
to be unique. A result of Paul Muhly, Jean Renault and Dana Williams \
establishes that the $C^{\ast }$-algebras of $G$ associated with two Haar
systems are strongly Morita equivalent (Theorem 2.8/p. 10 \cite{mrw}). If
the groupoid is transitive they have proved that if $G$ is transitive then
the\ $C^{\ast }$-algebra of $G$ is isomorphic to $C^{\ast }\left( H\right)
\otimes \mathcal{K}\left( L^{2}\left( \mu \right) \right) $, where $H$ is
the isotropy group $G_{u}^{u}$ at any unit $u\in G^{\left( 0\right) }$, $\mu 
$ is an essentially unique measure on $G^{\left( 0\right) }$, $C^{\ast
}\left( H\right) $ denotes the group $C^{\ast }$-algebra of $H$, and $%
\mathcal{K}\left( L^{2}\left( \mu \right) \right) $ denotes the compact
operators on $L^{2}\left( \mu \right) $ (Theorem 3.1/p. 16 \cite{mrw}).
Therefore the $C^{\ast }$-algebras of a \emph{transitive} groupoid $G$
associated with two Haar systems are $\ast $-isomorphic.

In \cite{rawm} Arlan Ramsay and Martin\ E.\ Walter have associated to a
locally compact groupoid $G$ a $C^{\ast }$- algebra denoted $M^{\ast }\left(
G,\nu \right) $. They have considered the universal representation $\omega $
of $C^{\ast }\left( G,\nu \right) $ -the usual $C^{\ast }$-algebra
associated to a Haar system $\nu =\left\{ \nu ^{u},\,u\in G^{\left( 0\right)
}\right\} $ (constructed as in \cite{re1}) . Since every cyclic
representation of $C^{\ast }\left( G,\nu \right) $ is the integrated form of
a representation of $G$, it follows that $\omega $ can be also regarded as a
representation of $\mathcal{B}_{c}\left( G\right) $, the space of compactly
supported Borel bounded function on $G$. Arlan Ramsay and Martin\ E.\ Walter
have used the notation $M^{\ast }\left( G,\nu \right) $ for the operator
norm closure of $\omega \left( \mathcal{B}_{c}\left( G\right) \right) $.
Since $\omega $ is an $\ast $-isomorphism on $C^{\ast }\left( G,\nu \right) $%
, we can regarded $C^{\ast }\left( G,\nu \right) $ as a subalgebra of $%
M^{\ast }\left( G,\nu \right) $.

We shall assume that the orbit space of the groupoid $G$ is proper and we
shall choose a Borel subset $F$ of $G^{\left( 0\right) }$ meeting each orbit
exactly once and such that $F\cap \left[ K\right] $ has a compact closure
for each compact subset $K$ of $G^{\left( 0\right) }$. For each appropriated
cross section $\sigma :G^{\left( 0\right) }\rightarrow G^{F}$ for $%
d_{F}:G^{F}\rightarrow G^{\left( 0\right) }$, $d_{F}\left( x\right) =d\left(
x\right) $, we shall construct a $C^{\ast }$-algebra $M_{\sigma }^{\ast
}\left( G,\nu \right) $ \ which can be viewed as a subalgebra of $M^{\ast
}\left( G,\nu \right) $. \ If $\nu _{1}=\left\{ \nu _{1}^{u},\,u\in
G^{\left( 0\right) }\right\} $ and $\nu _{2}=\left\{ \nu _{2}^{u},\,u\in
G^{\left( 0\right) }\right\} $ are two Haar systems on $G$, we shall prove
that the $C^{\ast }$-algebras $M_{\sigma }^{\ast }\left( G,\nu _{1}\right) $
and $M_{\sigma }^{\ast }\left( G,\nu _{2}\right) $ are $\ast $-isomorphic.

For a transitive (or more generally, a locally transitive) groupoid $G$ we
shall prove that the $C^{\ast }$-algebras $C^{\ast }\left( G,\nu \right) $, $%
M^{\ast }\left( G,\nu \right) $ and $M_{\sigma }^{\ast }\left( G,\nu \right) 
$ coincide.

For a principal proper groupoid $G$, we shall prove that 
\begin{equation*}
C^{\ast }\left( G,\nu \right) \subset M_{\sigma }^{\ast }\left( G,\nu
\right) \subset M^{\ast }\left( G,\nu \right) \text{.}
\end{equation*}

Let $\pi :G^{\left( 0\right) }\rightarrow G^{\left( 0\right) }/G$ be the
quotient map and let $\nu _{i}=\left\{ \varepsilon _{u}\times \mu
\,_{i}^{\pi \left( u\right) },\,u\in G^{\left( 0\right) }\right\} $, $%
\,i=1,2 $ be two Haar systems on the principal proper groupoid $G$. We shall
also prove that \ if the Hilbert bundles determined by the systems of
measures $\left\{ \mu _{i}^{\dot{u}}\right\} _{\dot{u}}$ have continuous
bases in the sense of Definition \ref{bcont}, then $\ast $-isomorphism
between $M_{\sigma }^{\ast }\left( G,\nu _{1}\right) $ and $M_{\sigma
}^{\ast }\left( G,\nu _{2}\right) $ can be restricted to a $\ast $%
-isomorphism between $C^{\ast }\left( G,\nu _{1}\right) $ and $C^{\ast
}\left( G,\nu _{2}\right) $.

For establishing notation, we include some definitions that can be found in
several places (e.g. \cite{re1}, \cite{mu}). A groupoid is a set $G$ endowed
with a product map 
\begin{equation*}
\left( x,y\right) \rightarrow xy\;\left[ :G^{\left( 2\right) }\rightarrow G%
\right]
\end{equation*}
where $G^{\left( 2\right) }$ is a subset of $G\times G$ called the set of
composable pairs, and an inverse map 
\begin{equation*}
x\rightarrow x^{-1}\;\left[ :G\rightarrow G\right]
\end{equation*}
such that the following conditions hold:

$\left( 1\right) \,$If $\left( x,y\right) \in G^{\left( 2\right) }$ and $%
\left( y,z\right) \in G^{\left( 2\right) }$, then $\left( xy,z\right) \in
G^{\left( 2\right) }$, $\left( x,yz\right) \in G^{\left( 2\right) }$ and $%
\left( xy\right) z=x\left( yz\right) $.

$\left( 2\right) $ $\left( x^{-1}\right) ^{-1}=x$ for all $x\in G$.

$\left( 3\right) $ For all $x\in G$, $\left( x,x^{-1}\right) \in G^{\left(
2\right) }$, and if $\left( z,x\right) \in G^{\left( 2\right) }$, then $%
\left( zx\right) x^{-1}=z$.

$\left( 4\right) $ For all $x\in G$, $\left( x^{-1},x\right) \in G^{\left(
2\right) }$, and if $\left( x,y\right) \in G^{\left( 2\right) }$, then $%
x^{-1}\left( xy\right) =y$.

The maps $r$ and $d$ on $G$, defined by the formulae $r\left( x\right)
=xx^{-1}$ and $d\left( x\right) =x^{-1}x$, are called the range and the
source maps. It follows easily from the definition that they have a common
image called the unit space of $G$, which is denoted $G^{\left( 0\right) }$.
Its elements are units in the sense that $xd\left( x\right) =r\left(
x\right) x=x$. Units will usually be denoted by letters as $u,\,v,\,w$ while
arbitrary elements will be denoted by $x,\,y,\,z$. It is useful to note that
a pair $\left( x,y\right) $ lies in $G^{\left( 2\right) }$ precisely when $%
d\left( x\right) =r\left( y\right) $, and that the cancellation laws hold
(e.g. $xy=xz$ iff $y=z$)$.$ The fibres of the range and the source maps are
denoted $G^{u}=r^{-1}\left( \left\{ u\right\} \right) $ and $%
G_{v}=d^{-1}\left( \left\{ v\right\} \right) $, respectively. More
generally, given the subsets $A$, $B$ $\subset G^{\left( 0\right) }$, we
define $G^{A}=r^{-1}\left( A\right) $, $G_{B}=d^{-1}\left( B\right) $ and $%
G_{B}^{A}=r^{-1}\left( A\right) \cap d^{-1}\left( B\right) $. The reduction
of $G$ to $A\subset G^{\left( 0\right) }$ is $G|A=G_{A}^{A}$. The relation $u%
\symbol{126}v$ iff $G_{v}^{u}\neq \phi $ is an equivalence relation on $%
G^{\left( 0\right) }$. Its equivalence classes are called orbits and the
orbit of a unit $u$ is denoted $\left[ u\right] $. A groupoid is called
transitive iff it has a single orbit. The quotient space for this
equivalence relation is called the orbit space of $G$ and denoted $G^{\left(
0\right) }/G$. We denote by $\pi :G^{\left( 0\right) }\rightarrow G^{\left(
0\right) }/G,\pi \left( u\right) =\dot{u}$ the quotient map. A subset of $%
G^{\left( 0\right) }$ is said saturated if it contains the orbits of its
elements. For any subset $A$ of $G^{\left( 0\right) }$, we denote by $\left[
A\right] $ the union of the orbits $\left[ u\right] $ for all $u\in A$.

A topological groupoid consists of a groupoid $G$ and a topology compatible
with the groupoid structure. This means that:

$\left( 1\right) \;x\rightarrow x^{-1}\;\left[ :G\rightarrow G\right] $ is
continuous.

$\left( 2\right) $\ $\left( x,y\right) \;\left[ :G^{\left( 2\right)
}\rightarrow G\right] $ is continuous where $G^{\left( 2\right) }$ has the
induced topology from $G\times G$.

We are exclusively concerned with topological groupoids which are second
countable, locally compact Hausdorff. It was shown in \cite{ra1} that
measured groupoids may be assume to have locally compact topologies, with no
loss in generality.

If $X$ is a locally compact space, $C_{c}\left( X\right) $ denotes the space
of complex-valuated continuous functions with compact support. The Borel
sets of a topological space are taken to be the $\sigma $-algebra generated
by the open sets. The space of compactly supported bounded Borel function on 
$X$ is denoted by $\mathcal{B}_{c}\left( X\right) $.

For a locally compact groupoid $G$, we denote by 
\begin{equation*}
G^{\prime }=\left( x\in G:\,r\left( x\right) =d\left( x\right) \right)
\end{equation*}
the isotropy group bundle of $G$. It is closed in $G$.

Let $G$ be a locally compact second countable groupoid equipped with a Haar
system, i.e. a family of positive Radon measures on $G$, $\left\{ \nu
^{u},\,u\in G^{\left( 0\right) }\right\} $, such that

$1)$ For all $u\in G^{\left( 0\right) }$, $supp(\nu ^{u})=\,G^{u}$.

$2)$ For all $f\in C_{c}\left( G\right) $, 
\begin{equation*}
u\rightarrow \int f\left( x\right) d\nu ^{u}\left( x\right) \,\;\left[
:G^{\left( 0\right) }\rightarrow \mathbf{C}\right]
\end{equation*}
is continuous.

$3)$ For all $f\in C_{c}\left( G\right) $ and all $x\in G$, 
\begin{equation*}
\int f\left( y\right) d\nu ^{r\left( x\right) }\left( y\right) \,=\,\int
f\left( xy\right) d\nu ^{d\left( x\right) }\left( y\right) \;
\end{equation*}

As a consequence of the existence of continuous Haar systems, $%
r,d:G\rightarrow G^{\left( 0\right) }$\ are open maps (\cite{we}).
Therefore, in this paper we shall always assume that $r:G\rightarrow
G^{\left( 0\right) }$\ is an open map

If $\mu $ is a measure on $G^{\left( 0\right) }$, then the measure $\nu
=\int \nu ^{u}d\mu \left( u\right) $, defined by 
\begin{equation*}
\int f\left( y\right) d\nu \left( y\right) =\int \left( \int f\left(
y\right) d\nu ^{u}\left( y\right) \right) d\mu \left( u\right) \text{, \ \ }%
f\geq 0\text{ Borel}
\end{equation*}
is called the measure on $G$ induced by $\mu $. The image of $\nu $ by the
inverse map $x\rightarrow x^{-1}$ is denoted $\nu ^{-1}$. $\mu $ is said
quasi-invariant if its induced measure $\nu $ is equivalent to its inverse $%
\nu ^{-1}$ \ A measure belongings to the class of \ a quasi-invariant
measure is also quasi-invariant. We say that the class is invariant.

If $\mu $ is a quasi-invariant measure on $G^{\left( 0\right) }$ and $\nu $
is the measure induced on $G$, then the Radon-Nikodym derivative $\Delta =%
\frac{d\nu }{d\nu ^{-1}}$ is called the modular function of $\mu $.

In order to define the $C^{\ast }$-algebra of a groupoid the space of
continuous functions with compact support on groupoid is made into a $\ast $%
-algebra and endowed with the smallest $C^{\ast }$-norm making its
representations continuous. For $f$, $g\in C_{c}\left( G\right) $ the
convolution is defined by:

\begin{equation*}
f\ast g\left( x\right) =\int f\left( xy\right) g\left( y^{-1}\right) d\nu
^{d\left( x\right) }\left( y\right)
\end{equation*}
and the involution by 
\begin{equation*}
f^{\ast }\left( x\right) =\overline{f\left( x^{-1}\right) }\text{.}
\end{equation*}
Under these operations, $C_{c}\left( G\right) $ becomes a topological $\ast $%
-algebra.

A representation of $C_{c}\left( G\right) $ is a $\ast $-homomorphism from $%
C_{c}\left( G\right) $ into $\mathcal{B}\left( H\right) $, for some Hilbert
space $H$, that is continuous with respect to the inductive limit topology
on $C_{c}\left( G\right) $ and the weak operator topology on $\mathcal{B}%
\left( H\right) $. The full $C^{\ast }$-algebra $C^{\ast }\left( G\right) $
is defined as the completion of the involutive algebra $C_{c}\left( G\right) 
$ with respect to the full $C^{\ast }$-norm 
\begin{equation*}
\left\| f\right\| =\sup \left\| L\left( f\right) \right\|
\end{equation*}
where $L$ runs over all non-degenerate representation of $C_{c}\left(
G\right) $ which are continuous for the inductive limit topology.

Every representation $\left( \mu ,G^{\left( 0\right) }\ast \mathcal{H}%
,L\right) $(see Definition 3.20/p.68 \cite{mu}) of $G$ can be integrated
into a representation, still denoted by $L$, of $C_{c}\left( G\right) $. The
relation between the two representation is: 
\begin{equation*}
\left\langle L\left( f\right) \xi _{1},\xi _{2}\right\rangle =\int f\left(
x\right) \left\langle L\left( x\right) \xi _{1}\left( d\left( x\right)
\right) ,\xi _{2}\left( r\left( x\right) \right) \right\rangle \Delta ^{-%
\frac{1}{2}}\left( x\right) d\nu ^{u}\left( x\right) d\mu \left( u\right)
\end{equation*}
where $f\in $ $C_{c}\left( G\right) $, $\xi _{1},\xi _{2}\in \int_{G^{\left(
0\right) }}^{\oplus }\mathcal{H}\left( u\right) d\mu \left( u\right) $.

Conversely, every non-degenerate $\ast $-representation of $C_{c}\left(
G\right) $ is obtained in this fashion (see \cite{re1} or \cite{mu}).

\section{The decomposition of a Haar system over the principal groupoid}

First we present some results on the structure of the Haar systems, as
developed by J. Renault in Section 1 of \cite{re2} and also by A. Ramsay and
M.E. Walter in Section 2 of \cite{rawm}.

In Section 1 of \cite{re2} Jean Renault constructs a Borel Haar system for $%
G^{\prime }$. One way to do this is to choose a function $F_{0}$ continuous
with conditionally support which is nonnegative and equal to $1$ at each $%
u\in G^{\left( 0\right) }.$ Then for each $u\in G^{\left( 0\right) }$ choose
a left Haar measure $\beta _{u}^{u}$ on $G_{u}^{u}$ so the integral of $%
F_{0} $ with respect to $\beta _{u}^{u}$ is $1.$

Renault defines $\beta _{v}^{u}=x\beta _{v}^{v}$ if $x\in G_{v}^{u}$ (where $%
x\beta _{v}^{v}\left( f\right) =\int f\left( xy\right) d\beta _{v}^{v}\left(
y\right) $ as usual). If $z$ is another element in $G_{v}^{u}$, then $%
x^{-1}z\in G_{v}^{v}$, and since $\beta _{v}^{v}$ is a left Haar measure on $%
G_{v}^{v}$, it follows that $\beta _{v}^{u}$ is independent of the choice of 
$x$. If $K$ is a compact subset of $G$, then $\sup\limits_{u,v}\beta
_{v}^{u}\left( K\right) <\infty $. Renault also defines a 1-cocycle $\delta $
on $G$ such that for every $u\in G^{\left( 0\right) }$, $\delta
|_{G_{u}^{u}} $ is the modular function for $\beta _{u}^{u}$. $\delta $ and $%
\delta ^{-1}=1/\delta $ are bounded on compact sets in $G$.

Let 
\begin{equation*}
R=\left( r,d\right) \left( G\right) =\left\{ \left( r\left( x\right)
,d\left( x\right) \right) ,\;x\in G\right\}
\end{equation*}
be the graph of the equivalence relation induced on $G^{\left( 0\right) }$.
This $R$ is the image of $G$ under the homomorphism $\left( r,d\right) $, so
it is a $\sigma $-compact groupoid. With this apparatus in place, Renault
describes a decomposition of the Haar system $\left\{ \nu ^{u},u\in
G^{\left( 0\right) }\right\} $ for $G$ over the equivalence relation $R$
(the principal groupoid associated to $G$). He proves that there is a unique
Borel Haar system $\alpha $ for $R$ \ with the property that 
\begin{equation*}
\nu ^{u}=\int \beta _{t}^{s}d\alpha ^{u}\left( s,t\right) \text{ \ for all }%
u\in G^{\left( 0\right) }\text{.}
\end{equation*}

In Section 2 \cite{rawm} A. Ramsay and M.E. Walter prove that 
\begin{equation*}
\sup\limits_{u}\alpha ^{u}\left( \left( r,d\right) \left( K\right) \right)
<\infty \text{, for all compact }K\subset G
\end{equation*}

For each $u\in G^{\left( 0\right) }$ the measure $\alpha ^{u}$ is
concentrated on $\left\{ u\right\} \times \left[ u\right] $. Therefore there
is a measure $\mu ^{u}$ concentrated on $\left[ u\right] $ such that $\alpha
^{u}=\varepsilon _{u}\times \mu ^{u}$, where $\varepsilon _{u}$ is the unit
point mass at $u$. Since $\left\{ \alpha ^{u},u\in G^{\left( 0\right)
}\right\} $ is a Haar system, we have $\mu ^{u}=\mu ^{v}$ for all $\left(
u,v\right) \in R$, and the function 
\begin{equation*}
u\rightarrow \int f\left( s\right) \mu ^{u}\left( s\right)
\end{equation*}
is Borel for all $f\geq 0$ Borel on $G^{\left( 0\right) }$. For each $u$ the
measure $\mu ^{u}$ is quasi-invariant (Section 2 \cite{rawm}). Therefore $%
\mu ^{u}$ is equivalent to $d_{\ast }\left( v^{u}\right) $ (Lemma 4.5/p. 277 
\cite{ra}).

If $\eta $ is a quasi-invariant measure for $\left\{ \nu ^{u},u\in G^{\left(
0\right) }\right\} $, then $\eta $ is a quasi-invariant measure for $\left\{
\alpha ^{u},u\in G^{\left( 0\right) }\right\} $. Also if $\Delta _{R}$ is
the modular function associated to $\left\{ \alpha ^{u},u\in G^{\left(
0\right) }\right\} $ and $\eta $, then $\Delta =\delta \Delta _{R}\circ
\left( r,d\right) $ can serve as the modular function associated to $\left\{
\nu ^{u},u\in G^{\left( 0\right) }\right\} $ and $\eta $.

Since $\mu ^{u}=\mu ^{v}$ for all $\left( u,v\right) \in R$, the system of
measures $\left\{ \mu ^{u}\right\} _{u}$ may be indexed on the elements of
the orbit space $G^{\left( 0\right) }/G$. \ 

\begin{definition}
We shall call the pair of the system of measures 
\begin{equation*}
\left( \left\{ \beta _{v}^{u}\right\} _{\left( u,v\right) \in R},\left\{ \mu
^{\dot{u}}\right\} _{\dot{u}\in G^{\left( 0\right) }/G}\right)
\end{equation*}
(described above) the decomposition of the Haar system $\left\{ \nu
^{u},u\in G^{\left( 0\right) }\right\} $ over the principal groupoid
associated to $G$. Also we shall call $\delta $ the $1$-cocycle associated
to the decomposition.
\end{definition}

\begin{remark}
Let us note that the system of measures $\left\{ \beta _{v}^{u}\right\} $
and the $1$-cocycle in the preceding definition do not depend on the Haar
system.
\end{remark}

\begin{lemma}
Let $G$ be a locally compact second countable groupoid with the bundle map $%
r|_{G^{\prime }}$ of $G^{\prime }$ open. Let $\left\{ \nu ^{u},u\in
G^{\left( 0\right) }\right\} $ be a Haar system on $G$ and $\left( \left\{
\beta _{v}^{u}\right\} ,\left\{ \mu ^{\dot{u}}\right\} \right) $ its
decomposition over the principal groupoid associated to $G$. Then for each $%
f\in C_{c}\left( G\right) $ the function 
\begin{equation*}
x\rightarrow \int f\left( y\right) d\beta _{d\left( x\right) }^{r\left(
x\right) }\left( y\right)
\end{equation*}
is continuous on $G$.
\end{lemma}

\begin{proof}
By Lemma 1.3/p. 6 \cite{re2}, for each $f\in C_{c}\left( G\right) $ the
function $u\rightarrow \int f\left( y\right) d\beta _{u}^{u}\left( y\right) $
is continuous$.$

\ Let $x\in G$ and $\left( x_{i}\right) _{i}$ be a sequence in $G$
converging to $x$. Let $f\in C_{c}\left( G\right) $ and let $g$ be a
continuous extension on $G$ of $y\rightarrow f\left( xy\right) \left[
:G^{d\left( x\right) }\rightarrow \mathbf{C}\right] $. Let $K$ be the
compact set 
\begin{equation*}
\left( \left\{ x,x_{i},i=1,2,..\right\} ^{-1}supp\left( f\right) \cup
supp\left( g\right) \right) \cap r^{-1}\left( \left\{ d\left( x\right)
,d\left( x_{i}\right) ,i=1,2,...\right\} \right) .
\end{equation*}
We have

$\left| \int f\left( y\right) d\beta _{d\left( x\right) }^{r\left( x\right)
}\left( y\right) -\int f\left( y\right) d\beta _{d\left( x_{i}\right)
}^{r\left( x_{i}\right) }\left( y\right) \right| $

$=\left| \int f\left( xy\right) d\beta _{d\left( x\right) }^{d\left(
x\right) }\left( y\right) -\int f\left( x_{i}y\right) d\beta _{d\left(
x_{i}\right) }^{d\left( x_{i}\right) }\left( y\right) \right| $

$=\left| \int g\left( y\right) d\beta _{d\left( x\right) }^{d\left( x\right)
}\left( y\right) -\int f\left( x_{i}y\right) d\beta _{d\left( x_{i}\right)
}^{d\left( x_{i}\right) }\left( y\right) \right| $

$\leq \left| \int g\left( y\right) d\beta _{d\left( x\right) }^{d\left(
x\right) }\left( y\right) -\int g\left( y\right) d\beta _{d\left(
x_{i}\right) }^{d\left( x_{i}\right) }\left( y\right) \right| +$

$\ \ \ \ \ \ \ \ \ \ \ \ \ \ \ \ \ \ \ \ \ \ +\left| \int g\left( y\right)
d\beta _{d\left( x_{i}\right) }^{d\left( x_{i}\right) }\left( y\right) -\int
f\left( x_{i}y\right) d\beta _{d\left( x_{i}\right) }^{d\left( x_{i}\right)
}\left( y\right) \right| $

$\leq \left| \int g\left( y\right) d\beta _{d\left( x\right) }^{d\left(
x\right) }\left( y\right) -\int g\left( y\right) d\beta _{d\left(
x_{i}\right) }^{d\left( x_{i}\right) }\left( y\right) \right| +$

$\ \ \ \ \ \ \ \ \ \ \ \ \ \ \ \ \ \ \;\;\;+\sup_{y\in G_{d\left(
x_{i}\right) }^{d\left( x_{i}\right) }}\left| g\left( y\right) -f\left(
x_{i}y\right) \right| \beta _{d\left( x_{i}\right) }^{d\left( x_{i}\right)
}\left( K\right) $

A compactness argument shows that $\sup_{y\in G_{d\left( x_{i}\right)
}^{d\left( x_{i}\right) }}\left| g\left( y\right) -f\left( x_{i}y\right)
\right| $ converges to $0$. Also $\left| \int g\left( y\right) d\beta
_{d\left( x\right) }^{d\left( x\right) }\left( y\right) -\int g\left(
y\right) d\beta _{d\left( x_{i}\right) }^{d\left( x_{i}\right) }\left(
y\right) \right| $ converges to $0$ because the function $u\rightarrow \int
f\left( y\right) d\beta _{u}^{u}\left( y\right) $ is continuous$.$ Hence 
\begin{equation*}
\left| \int f\left( y\right) d\beta _{d\left( x\right) }^{r\left( x\right)
}\left( y\right) -\int f\left( y\right) d\beta _{d\left( x_{i}\right)
}^{r\left( x_{i}\right) }\left( y\right) \right|
\end{equation*}
converges to $0.$
\end{proof}

\begin{definition}
A locally compact groupoid $G$ is proper if the map $\left( r,d\right)
:G\rightarrow G^{\left( 0\right) }\times G^{\left( 0\right) }$ is proper
(i.e. the inverse map of each compact subset of $G^{\left( 0\right) }\times
G^{\left( 0\right) }$ is compact).(Definition 2.1.9/p. 37 \cite{ar}).
\end{definition}

Throughout this paper we shall assume that $G$ is a second countable locally
compact groupoid for which the orbit space is Hausdorff and the map 
\begin{equation*}
\left( r,d\right) :G\rightarrow R,\,\left( r,d\right) \left( x\right)
=\left( r\left( x\right) ,d\left( x\right) \right)
\end{equation*}
is open, where $R$ is endowed with the product topology induced from $%
G^{\left( 0\right) }\times G^{\left( 0\right) }$. Therefore $R$ will be a
locally compact groupoid. The fact that $R$ is a closed subset of $G^{\left(
0\right) }\times G^{\left( 0\right) }$ and that it is endowed with the
product topology is equivalent to the fact $R$ is a proper groupoid.

Throughout this paper by a groupoid with proper orbit space we shall mean a
groupoid $G$ for which the orbit space is Hausdorff and the map 
\begin{equation*}
\left( r,d\right) :G\rightarrow R,\,\left( r,d\right) \left( x\right)
=\left( r\left( x\right) ,d\left( x\right) \right)
\end{equation*}
is open, where $R$ is endowed with the product topology induced from $%
G^{\left( 0\right) }\times G^{\left( 0\right) }$.

\begin{proposition}
\label{pgr}Let $G$ be a second countable locally compact groupoid with
proper orbit space. Let $\left\{ \nu ^{u},u\in G^{\left( 0\right) }\right\} $
be a Haar system on $G$ and $\left( \left\{ \beta _{v}^{u}\right\} ,\left\{
\mu ^{\dot{u}}\right\} \right) $ its decomposition over the principal
groupoid associated to $G$. Then .for each $g\in C_{c}\left( G^{\left(
0\right) }\right) $, the map 
\begin{equation*}
u\rightarrow \int g\left( v\right) d\mu ^{\pi \left( u\right) }\left(
v\right)
\end{equation*}
is continuous.
\end{proposition}

\begin{proof}
Let $g\in C_{c}\left( G^{\left( 0\right) }\right) $ and $u_{0}\in G^{\left(
0\right) }$. Let $K_{1}$ be a compact neighborhood of $u_{0}$ and $K_{2}$ be
the support of $g$. Since \ $G$ is locally compact and $\left( r,d\right) $
is open from $G$ to $\left( r,d\right) \left( G\right) $, there is a compact
subset $K$ of $G$ such that $\left( r,d\right) \left( K\right) $ contains $%
\left( K_{1}\times K_{2}\right) \cap \left( r,d\right) \left( G\right) $.
Let $F_{1}\in C_{c}\left( G\right) $ be a nonnegative function equal to $1$
on a compact neighborhood $U$ of $K$. Let $F_{2}$ $\in C_{c}\left( G\right) $
be a function which extends to $G$ the function $x\rightarrow F_{1}\left(
x\right) /\int F_{1}\left( y\right) d\beta _{d\left( x\right) }^{r\left(
x\right) }\left( y\right) $, $x\in U$. We have $\int F_{2}\left( y\right)
d\beta _{v}^{u}\left( y\right) =1$ for all $\left( u,v\right) \in \left(
r,d\right) \left( K\right) $. Since for all $u\in K_{1}$, 
\begin{eqnarray*}
\int g\left( v\right) d\mu ^{\pi \left( u\right) }\left( v\right) &=&\int
g\left( v\right) \int F_{2}\left( y\right) d\beta _{v}^{u}\left( y\right)
d\mu ^{\pi \left( u\right) }\left( v\right) \\
&=&\int g\left( d\left( y\right) \right) F_{2}\left( y\right) d\nu
^{u}\left( y\right) \text{,}
\end{eqnarray*}
it follows that $u\rightarrow \int g\left( v\right) d\mu ^{\pi \left(
u\right) }\left( v\right) $ is continuous at $u_{0}.$
\end{proof}

\begin{remark}
Let $G$ be a locally compact second countable groupoid with proper orbit
space. Let $\left\{ \nu ^{u},u\in G^{\left( 0\right) }\right\} $be a Haar
system on $G$ and $\left( \left\{ \beta _{v}^{u}\right\} ,\left\{ \mu ^{\dot{%
u}}\right\} \right) $ be its decomposition over the associated principal
groupoid. If $\mu $ is a quasi-invariant probability measure for the Haar
system, then $\mu _{1}=\int \mu ^{\pi \left( u\right) }d\mu \left( u\right) $
is a Radon measure which is equivalent to $\mu .$ Indeed, let $f\geq 0$
Borel on $G^{\left( 0\right) }$ such that $\mu \left( f\right) =0$. Since $%
\mu $ is quasi-invariant, it follows that for $\mu \,a.a.\;u$ $\nu
^{u}\left( f\circ d\right) =0$, and since $\mu ^{\pi \left( u\right) }$ is
equivalent to $d_{\ast }\left( v^{u}\right) $, it results $\mu ^{\pi \left(
u\right) }\left( f\right) =0$ for $\mu $ $a.a.\;u$. Conversely, if $\mu
_{1}\left( f\right) =0$, then $\mu ^{\pi \left( u\right) }\left( f\right) =0$
for $\mu $ $a.a.\;u$, and therefore $\nu ^{u}\left( f\circ d\right) =0$.
Thus the quasi-invariance of $\mu $ implies $\mu \left( f\right) =0$. Thus
each Radon quasi-invariant measure is equivalent to a Radon measure of the
form $\int \mu ^{\dot{u}}d\tilde{\mu}\left( \dot{u}\right) $, where $\tilde{%
\mu}$ is a probability measure on the orbit space $G/G^{\left( 0\right) }.$
\end{remark}

\section{A $C^{\ast }$-algebra associated to a locally compact groupoid with
proper orbit space}

Let $G$ be a locally compact second countable groupoid with proper orbit
space. Let 
\begin{equation*}
\pi :G^{\left( 0\right) }\rightarrow G^{\left( 0\right) }/G
\end{equation*}
be the quotient map. Since the quotient space is proper, $G^{\left( 0\right)
}/G$ is Hausdorff. Let us assume that the range map $r$ is open. As a
consequence, the map $\pi $ is open. Applying Lemma 1.1 \cite{ma} to the
locally compact second countable spaces $G^{\left( 0\right) }$ and $%
G^{\left( 0\right) }/G$ and to the continuous open surjection $\pi
:G^{\left( 0\right) }\rightarrow G^{\left( 0\right) }/G$, it follows that
there is a Borel set $F$ in $G^{\left( 0\right) }$ such that:

\begin{enumerate}
\item  $F$ contains exactly one element in each orbit $\left[ u\right] =\pi
^{-1}\left( \pi \left( u\right) \right) $.

\item  For each compact subset $K$ of $G^{\left( 0\right) }$, $F\cap \left[ K%
\right] =F\cap \pi ^{-1}\left( \pi \left( K\right) \right) $ has a compact
closure.
\end{enumerate}

For each unit $u$ let us define $e\left( u\right) =F\cap \left[ u\right] $ ($%
e\left( u\right) $ is the unique element in the orbit of $u$ contained in $F$%
). For each Borel subset $B$ of $G^{\left( 0\right) }$, $\pi $ is continuous
and one-to-one on $B\cap F$ and hence $\pi \left( B\cap F\right) $ is Borel
in $G^{\left( 0\right) }/G$. Therefore the map $e:G^{\left( 0\right)
}\rightarrow G^{\left( 0\right) }$ is Borel (for each Borel subset $B$ of $%
G^{\left( 0\right) }$, $e^{-1}\left( B\right) =\left[ B\cap F\right] =\pi
^{-1}\left( \pi \left( B\cap F\right) \right) $ is Borel in $G^{\left(
0\right) }$). Also for each compact subset $K$ of $G^{\left( 0\right) }$, $%
e\left( K\right) $ has a compact closure because $e\left( K\right) \subset
F\cap \left[ K\right] $.

Since the orbit space $G^{\left( 0\right) }/G$ is proper the map 
\begin{equation*}
\left( r,d\right) :G\rightarrow R,\,\left( r,d\right) \left( x\right)
=\left( r\left( x\right) ,d\left( x\right) \right)
\end{equation*}
is open and $R$ is closed in $G^{\left( 0\right) }\times G^{\left( 0\right)
} $. Applying Lemma 1.1 \cite{ma} to the locally compact second countable
spaces $G$ and $R$ and to the continuous open surjection $\left( r,d\right)
:G\rightarrow R$, it follows that there is a\textit{\ regular cross section} 
$\sigma _{0}:R\rightarrow G$. This means that $\sigma _{0}$ is Borel, $%
\left( r,d\right) \left( \sigma _{0}\left( u,v\right) \right) =\left(
u,v\right) $ for all $\left( u,v\right) \in R$, and $\sigma _{0}\left(
K\right) $ is relatively compact in $G$ for each compact subset $K$ of $R$.

Let us define $\sigma :G^{\left( 0\right) }\rightarrow G^{F}$ by $\sigma
\left( u\right) =\sigma _{0}\left( e\left( u\right) ,u\right) $ for all $u$.
It is easy to note that $\sigma $ is a\textit{\ }cross section for $%
d:G^{F}\rightarrow G^{\left( 0\right) }$ and $\sigma \left( K\right) $ is
relatively compact in $G$ for all compact $K\subset G^{\left( 0\right) }$.

Replacing $\sigma $ by 
\begin{equation*}
v\rightarrow \sigma \left( e\left( v\right) \right) ^{-1}\sigma \left(
v\right)
\end{equation*}
we may assume that $\sigma \left( e\left( v\right) \right) =e\left( v\right) 
$ for all $v$. Let us define $q:G\rightarrow G_{F}^{F}$ by 
\begin{equation*}
q\left( x\right) =\sigma \left( r\left( x\right) \right) x\sigma \left(
d\left( x\right) \right) ^{-1},\,x\in G\text{.}
\end{equation*}

Let $\nu =\left\{ \nu ^{u}:\,u\in G^{\left( 0\right) }\right\} $ be a Haar
system on $G$ and let $\left( \left\{ \beta _{v}^{u}\right\} ,\left\{ \mu ^{%
\dot{u}}\right\} \right) $ be its decompositions over the principal groupoid
. Let $\delta $ be the $1$-cocycle associated to the decomposition.

Let us denote by $\mathcal{B}_{\sigma }\left( G\right) $ the linear span of
the functions of the form 
\begin{equation*}
x\rightarrow g_{1}\left( r\left( x\right) \right) g\left( q\left( x\right)
\right) g_{2}\left( d\left( x\right) \right)
\end{equation*}
where $g_{1},g_{2}$ are compactly supported Borel bounded function on $%
G^{\left( 0\right) }$ and $g$ is a Borel bounded function on $G_{F}^{F}$
such that if $S$ is the support of $g$, then the closure of $\sigma \left(
K_{1}\right) ^{-1}S\sigma \left( K_{2}\right) $ is compact in $G$ for all
compact subsets $K_{1}$, $K_{2}$ of $G^{\left( 0\right) }$. \ $\mathcal{B}%
_{\sigma }\left( G\right) $ is a subspace of $\mathcal{B}_{c}\left( G\right) 
$, the space of compactly supported Borel bounded function on $G$.

If $f_{1}$, $f_{2}\in \mathcal{B}_{\sigma }\left( G\right) $ are defined by 
\begin{eqnarray*}
f_{1}\left( x\right) &=&g_{1}\left( r\left( x\right) \right) g\left( q\left(
x\right) \right) g_{2}\left( d\left( x\right) \right) \\
f_{2}\left( x\right) &=&h_{1}\left( r\left( x\right) \right) h\left( q\left(
x\right) \right) h_{2}\left( d\left( x\right) \right)
\end{eqnarray*}
then 
\begin{equation*}
f_{1}\ast f_{2}\left( x\right) =g\ast h\left( q\left( x\right) \right)
g_{1}\left( r\left( x\right) \right) h_{2}\left( d\left( x\right) \right)
\left\langle g_{2},\overline{h_{1}}\right\rangle _{\pi \left( r\left(
x\right) \right) }
\end{equation*}
\begin{equation*}
f_{1}^{\ast }\left( x\right) =\overline{g_{2}\left( r\left( x\right) \right) 
}\overline{g\left( q\left( x\right) ^{-1}\right) }\overline{g_{1}\left(
d\left( x\right) \right) }
\end{equation*}
Thus $\mathcal{B}_{\sigma }\left( G\right) $ is closed under convolution and
involution.

Let $\omega $ be the universal representation of $C^{\ast }\left( G,\nu
\right) $ the usual $C^{\ast }$-algebra associated to a Haar system $\nu
=\left\{ \nu ^{u},\,u\in G^{\left( 0\right) }\right\} $ (constructed as in 
\cite{re1}) . Since every cyclic representation of $C^{\ast }\left( G,\nu
\right) $ is the integrated form of a representation of $G$, it follows that 
$\omega $ can be also regarded as a representation of $\mathcal{B}_{c}\left(
G\right) $, the space of compactly supported Borel bounded function on $G$.
Arlan Ramsay and Martin\ E.\ Walter have used the notation $M^{\ast }\left(
G,\nu \right) $ for the operator norm closure of $\omega \left( \mathcal{B}%
_{c}\left( G\right) \right) $. Since $\omega $ is an $\ast $-isomorphism on $%
C^{\ast }\left( G,\nu \right) $, we can regarded $C^{\ast }\left( G,\nu
\right) $ as a subalgebra of $M^{\ast }\left( G,\nu \right) $.

\begin{definition}
We denote by $M_{\sigma }^{\ast }\left( G,\nu \right) $ the operator norm
closure of $\omega \left( \mathcal{B}_{\sigma }\left( G\right) \right) $.
\end{definition}

\begin{lemma}
\label{uni}Let $\left\{ \mu _{1}^{\dot{u}}\right\} _{\dot{u}}$ and $\left\{
\mu _{2}^{\dot{u}}\right\} _{\dot{u}}$ be two systems of measures on $%
G^{\left( 0\right) }$ satisfying:

\begin{enumerate}
\item  $supp\left( \mu _{i}^{\dot{u}}\right) =\left[ u\right] $ for all $%
\dot{u}$, $i=1,2$

\item  For all compactly supported Borel bounded function $f$ on $G^{\left(
0\right) }$ the function 
\begin{equation*}
u\rightarrow \int f\left( v\right) \mu _{i}^{\pi \left( u\right) }\left(
v\right)
\end{equation*}
is bounded and Borel.
\end{enumerate}

Then there is a family $\left\{ U_{\dot{u}}\right\} _{\dot{u}}$ of unitary
operators with the following properties:

\begin{enumerate}
\item  $U_{\dot{u}}:L^{2}\left( \mu _{1}^{\dot{u}}\right) \rightarrow
L^{2}\left( \mu _{2}^{\dot{u}}\right) $ is a unitary operator for each $\dot{%
u}\in G^{\left( 0\right) }/G$.

\item  For all Borel bounded function $f$ on $G^{\left( 0\right) }$, 
\begin{equation*}
u\rightarrow U_{\pi \left( u\right) }\left( f\right)
\end{equation*}
is a bounded Borel function with compact support.

\item  For all Borel bounded function $f$ on $G^{\left( 0\right) }$, 
\begin{equation*}
U_{\pi \left( u\right) }\left( \overline{f}\right) =\overline{U_{\pi \left(
u\right) }\left( f\right) }
\end{equation*}
\end{enumerate}
\end{lemma}

\begin{proof}
Using the same argument as in \ \cite{ra1} (p. 323) we can construct a
sequence $f_{1},f_{2},...$ of real valued Borel bounded function on $%
G^{\left( 0\right) }$ such that $\dim \left( L^{2}\left( \mu _{1}^{\dot{u}%
}\right) \right) =\infty $ if and only if $\left\| f_{n}\right\| _{2}$ for $%
n=1,2,...$ and then $\left\{ f_{1},f_{2},...\right\} $ gives an ortonormate
basis of $L^{2}\left( \mu _{1}^{\dot{u}}\right) $, while $\dim \left(
L^{2}\left( \mu _{1}^{\dot{u}}\right) \right) =k<\infty $ if and only if $%
\left\| f_{n}\right\| _{2}=1$ for $n\leq k$, and $\left\| f_{n}\right\|
_{2}=0$ for $n>k$ \ and then $\left\{ f_{1},f_{2},...f_{k}\right\} $ gives
an ortonormate basis of $L^{2}\left( \mu _{1}^{\dot{u}}\right) $. \ Let \ $%
g_{1},g_{2},...$ be a sequence with the same properties as $f_{1}.f_{2},...$
\ corresponding to $\left\{ \mu _{2}^{\dot{u}}\right\} _{\dot{u}}$. Let us
define $U_{\dot{u}}:L^{2}\left( \mu _{1}^{\dot{u}}\right) \rightarrow
L^{2}\left( \mu _{2}^{\dot{u}}\right) $ by 
\begin{equation*}
U_{\dot{u}}\left( f_{n}\right) =g_{n}\text{ for all }n
\end{equation*}
Then the family $\left\{ U_{\dot{u}}\right\} _{\dot{u}}$ has the required
properties.
\end{proof}

\begin{theorem}
\label{iso}Let $G$ be a locally compact second countable groupoid with
proper orbit space. Let $\left\{ \nu _{i}^{u},u\in G^{\left( 0\right)
}\right\} $, $i=1,2$ be two Haar systems on $G$. Let $F$ be a Borel subset
of $G^{\left( 0\right) }$ containing only one element $e\left( u\right) $ in
each orbit $\left[ u\right] $. Let $\sigma :G^{\left( 0\right) }\rightarrow
G^{F}$ be a\ cross section for $d:G^{F}\rightarrow G^{\left( 0\right) }$
with $\sigma \left( e\left( v\right) \right) =e\left( v\right) $ for all $%
v\in G^{\left( 0\right) }$ and $\sigma \left( K\right) $ relatively compact
in $G$ for all compact sets $K\subset G^{\left( 0\right) }$ . Then the $%
C^{\ast }$-algebras $M_{\sigma }^{\ast }\left( G,\nu _{1}\right) $ and $%
M_{\sigma }^{\ast }\left( G,\nu _{2}\right) $ are $\ast $-isomorphic.
\end{theorem}

\begin{proof}
Let $\left( \left\{ \beta _{v}^{u}\right\} ,\left\{ \mu _{i}^{\dot{u}%
}\right\} \right) $ be the decompositions\ of the Haar systems over the
principal groupoid . Let $\delta $ be the $1$-cocycle associated to the
decompositions, $i=1,2$.

We shall denote by $\left\langle \cdot ,\,\cdot \right\rangle _{i,\dot{u}}$
the inner product of $\left( L^{2}\left( G^{\left( 0\right) },\,\delta
\left( \sigma \left( \cdot \right) \right) \mu _{i}^{\dot{u}}\right) \right) 
$, $i=1,2$.

Let us define $q:G\rightarrow G_{F}^{F}$ by 
\begin{equation*}
q\left( x\right) =\sigma \left( r\left( x\right) \right) x\sigma \left(
d\left( x\right) \right) ^{-1},\,x\in G\text{.}
\end{equation*}

We shall define a $\ast $-homomorphism $\Phi $ from $\mathcal{B}_{\sigma
}\left( G\right) $ to $\mathcal{B}_{\sigma }\left( G\right) $. It suffices
to define $\Phi $ on the set of function on $G$ of the form 
\begin{equation*}
x\rightarrow g_{1}\left( r\left( x\right) \right) g\left( q\left( x\right)
\right) g_{2}\left( d\left( x\right) \right)
\end{equation*}
Let $\left\{ U_{\dot{u}}\right\} _{\dot{u}}$ be the family of unitary
operators with the properties stated in Lemma \ref{uni}, associated to the
systems of measures $\left\{ \delta \left( \sigma \left( \cdot \right)
\right) \mu _{i}^{\dot{u}}\right\} _{\dot{u}}$, $i=1,2$.

Let us define $\Phi $ by 
\begin{equation*}
\Phi \left( f\right) =\left( x\rightarrow U_{\pi \left( r\left( x\right)
\right) }\left( g_{1}\right) \left( r\left( x\right) \right) g\left( q\left(
x\right) \right) U_{\pi \left( d\left( x\right) \right) }\left( g_{2}\right)
\left( d\left( x\right) \right) \right)
\end{equation*}
where $f$ is defined by 
\begin{equation*}
f\left( x\right) =g_{1}\left( r\left( x\right) \right) g\left( q\left(
x\right) \right) g_{2}\left( d\left( x\right) \right)
\end{equation*}

If $f_{1}$ and $f_{2}$ are defined by 
\begin{eqnarray*}
f_{1}\left( x\right) &=&g_{1}\left( r\left( x\right) \right) g\left( q\left(
x\right) \right) g_{2}\left( d\left( x\right) \right) \\
f_{2}\left( x\right) &=&h_{1}\left( r\left( x\right) \right) h\left( q\left(
x\right) \right) h_{2}\left( d\left( x\right) \right)
\end{eqnarray*}
then 
\begin{equation*}
f_{1}\ast f_{2}\left( x\right) =g\ast h\left( q\left( x\right) \right)
g_{1}\left( r\left( x\right) \right) h_{2}\left( d\left( x\right) \right)
\left\langle g_{2},\overline{h_{1}}\right\rangle _{1,\pi \left( r\left(
x\right) \right) }
\end{equation*}
and consequently 
\begin{eqnarray*}
&&\Phi \left( f_{1}\ast f_{2}\right) \\
&=&g\ast h\left( q\left( x\right) \right) U_{\pi \left( r\left( x\right)
\right) }\left( g_{1}\right) \left( r\left( x\right) \right) U_{\pi \left(
r\left( x\right) \right) }\left( h_{2}\right) \left( d\left( x\right)
\right) \left\langle g_{2},\overline{h_{1}}\right\rangle _{1,\pi \left(
r\left( x\right) \right) } \\
&=&\Phi \left( f_{1}\right) \ast \Phi \left( f_{2}\right) \text{.}
\end{eqnarray*}

Let $\tilde{\eta}$ be a probability measure on $G^{\left( 0\right) }/G$ and $%
\eta _{i}=\int \mu _{i}^{\dot{u}}d\tilde{\eta}\left( \dot{u}\right) ,i=1,2$.
Let $L_{1}$ be the integrated form of a representation $\left( L,\mathcal{H}%
\ast G^{\left( 0\right) },\eta _{1}\right) $ and $L_{2}$ be the integrated
form of $\left( L,\mathcal{H}\ast G^{\left( 0\right) },\eta _{2}\right) $.
Let $B$ be the Borel function defined by: 
\begin{equation*}
B\left( u\right) =L\left( \sigma \left( u\right) \right)
\end{equation*}
and $W:\int_{G^{\left( 0\right) }}^{\oplus }\mathcal{H}\left( u\right) d\eta
_{1}\left( u\right) \rightarrow \int_{G^{\left( 0\right) }}^{\oplus }$ $%
\mathcal{H}\left( e\left( u\right) \right) d\eta _{1}\left( u\right) $ be
defined by 
\begin{equation*}
W\left( \zeta \right) =\left( u\rightarrow B\left( u\right) \left( \zeta
\left( u\right) \right) \right)
\end{equation*}

Since every element of $L^{2}\left( G^{\left( 0\right) },\delta \left(
\sigma \left( \cdot \right) \right) \mu _{1}^{\dot{w}},\mathcal{H}\left(
e\left( w\right) \right) \right) $ is a limit of linear combinations of
elements $u\rightarrow a\left( u\right) \xi $ with $a\in L^{2}\left(
G^{\left( 0\right) },\delta \left( \sigma \left( \cdot \right) \right) \mu
_{1}^{\dot{w}}\right) $ and $\xi \in \mathcal{H}\left( e\left( w\right)
\right) $, we can define a unitary operator 
\begin{equation*}
V_{\dot{w}}:L^{2}\left( G^{\left( 0\right) },\delta \left( \sigma \left(
\cdot \right) \right) \mu _{1}^{\dot{w}},\mathcal{H}\left( e\left( w\right)
\right) \right) \rightarrow L^{2}\left( G^{\left( 0\right) },\delta \left(
\sigma \left( \cdot \right) \right) \mu _{2}^{\dot{w}},\mathcal{H}\left(
e\left( w\right) \right) \right)
\end{equation*}
by 
\begin{equation*}
V_{\dot{w}}\left( u\rightarrow a\left( u\right) \xi \right) =U_{\dot{w}%
}\left( a\right) \xi
\end{equation*}

Let $V:\int_{G^{\left( 0\right) }}^{\oplus }\mathcal{H}\left( e\left(
u\right) \right) d\eta _{1}\left( u\right) \rightarrow \int_{G^{\left(
0\right) }}^{\oplus }$ $\mathcal{H}\left( e\left( u\right) \right) d\eta
_{2}\left( u\right) $ be defined by 
\begin{equation*}
V\left( \zeta \right) =\left( u\rightarrow V_{\dot{u}}\left( \zeta \left(
u\right) \right) \right)
\end{equation*}

If $\zeta _{1},\zeta _{2}\in \int_{G^{\left( 0\right) }}^{\oplus }$ $%
\mathcal{H}\left( e\left( u\right) \right) d\eta _{1}\left( u\right) $ and $%
f $ is of the form 
\begin{equation*}
f\left( x\right) =g_{1}\left( r\left( x\right) \right) g\left( q\left(
x\right) \right) g_{2}\left( d\left( x\right) \right) \text{,}
\end{equation*}
we have 
\begin{eqnarray*}
&&\left\langle WL_{1}\left( f\right) W^{\ast }\zeta _{1},\zeta
_{2}\right\rangle \\
&=&\int \int g\left( x\right) \delta \left( x\right) ^{\frac{-1}{2}%
}\left\langle L\left( x\right) A_{1}\left( \dot{w}\right) ,B_{1}\left( \dot{w%
}\right) \right\rangle d\beta _{e\left( w\right) }^{e\left( w\right) }\left(
x\right) d\tilde{\eta}\left( \dot{w}\right)
\end{eqnarray*}
where 
\begin{eqnarray*}
A_{1}\left( \dot{w}\right) &=&\int g_{2}\left( v\right) \zeta _{1}\left(
v\right) \delta \left( \sigma \left( v\right) \right) ^{\frac{1}{2}}d\mu
_{1}^{\dot{w}}\left( v\right) \\
B_{1}\left( \dot{w}\right) &=&\int g_{1}\left( u\right) \zeta _{2}\left(
u\right) \delta \left( \sigma \left( u\right) \right) ^{\frac{1}{2}}d\mu
_{1}^{\dot{w}}\left( u\right)
\end{eqnarray*}

Moreover, if $f$ is of the form $f\left( x\right) =g_{1}\left( r\left(
x\right) \right) g\left( q(x\right) g_{2}\left( d\left( x\right) \right) $
and $\zeta _{1},\zeta _{2}\in \int_{G^{\left( 0\right) }}^{\oplus }$ $%
\mathcal{H}\left( e\left( u\right) \right) d\eta _{2}\left( u\right) $, then 
\begin{eqnarray*}
&&\left\langle VWL_{1}\left( f\right) W^{\ast }V^{\ast }\zeta _{1},\zeta
_{2}\right\rangle \\
&=&\int \int g\left( x\right) \delta \left( x\right) ^{\frac{-1}{2}%
}\left\langle L\left( x\right) A_{2}\left( \dot{w}\right) ,B_{2}\left( \dot{w%
}\right) \right\rangle d\beta _{e\left( u\right) }^{e\left( u\right) }\left(
x\right) d\tilde{\eta}\left( \dot{w}\right) \\
&=&\left\langle WL_{2}\left( \Phi \left( f\right) \right) W^{\ast }\zeta
_{1},\zeta _{2}\right\rangle
\end{eqnarray*}
where 
\begin{eqnarray*}
A_{2}\left( \dot{w}\right) &=&\int g_{2}\left( v\right) V^{\ast }\zeta
_{1}\left( v\right) \delta \left( \sigma \left( v\right) \right) ^{\frac{1}{2%
}}d\mu _{1}^{\dot{w}}\left( v\right) \\
&=&\int U_{\dot{v}}\left( g_{2}\right) \left( v\right) \zeta _{1}\left(
v\right) \delta \left( \sigma \left( v\right) \right) ^{\frac{1}{2}}d\mu
_{2}^{\dot{w}}\left( v\right) \\
B_{2}\left( \dot{w}\right) &=&\int g_{1}\left( v\right) V^{\ast }\zeta
_{2}\left( v\right) \delta \left( \sigma \left( v\right) \right) ^{\frac{1}{2%
}}d\mu _{1}^{\dot{w}}\left( v\right) \\
&=&\int U_{\dot{u}}\left( g_{1}\right) \left( u\right) \zeta _{2}\left(
u\right) \delta \left( \sigma \left( u\right) \right) ^{\frac{1}{2}}d\mu
_{2}^{\dot{w}}\left( u\right)
\end{eqnarray*}

Therefore $\left\| L_{1}\left( f\right) \right\| =\left\| L_{2}\left( \Phi
\left( f\right) \right) \right\| $. Consequently we can extend $\Phi $ to a $%
\ast $-homomorphism between the $M_{\sigma }^{\ast }\left( G,\nu _{1}\right) 
$ and $M_{\sigma }^{\ast }\left( G,\nu _{2}\right) $. It is not hard to see
that $\Phi $ is in fact a $\ast $-isomorphism: 
\begin{equation*}
\Phi ^{-1}\left( f\right) =\left( x\rightarrow U_{\pi \left( r\left(
x\right) \right) }^{\ast }\left( g_{1}\right) \left( r\left( x\right)
\right) g\left( q\left( x\right) \right) U_{\pi \left( d\left( x\right)
\right) }^{\ast }\left( g_{2}\right) \left( d\left( x\right) \right) \right)
\end{equation*}
for each $f$ of the form 
\begin{equation*}
f\left( x\right) =g_{1}\left( r\left( x\right) \right) g\left( q\left(
x\right) \right) g_{2}\left( d\left( x\right) \right) .
\end{equation*}
\end{proof}

\bigskip

\section{The case of locally transitive groupoids}

A locally compact transitive groupoid $G$ is a groupoid \ for which all
orbits $\left[ u\right] $ are open in $G^{\left( 0\right) }$. We shall prove
that if $G$ is a locally compact second countable locally transitive
groupoid endowed with a Haar system $\nu $, then

\begin{equation*}
C^{\ast }\left( G,\nu \right) =M^{\ast }\left( G,\nu \right) =M_{\sigma
}^{\ast }\left( G,\nu \right)
\end{equation*}
for any regular cross section $\sigma .$

\begin{notation}
Let $\left\{ \nu ^{u},\,u\in G^{\left( 0\right) }\right\} $ be a fixed Haar
system on $G$. Let $\mu $ be a quasi-invariant measure, $\Delta $ its
modular function, $\nu _{1}$ be the measure induced by $\mu $ on $G$ and $%
\nu _{0}=\Delta ^{-\frac{1}{2}}\nu _{1}$. Let us denote by $II_{\mu }\left(
G\right) $ the set 
\begin{equation*}
\left\{ f\in L^{1}\left( G,\nu _{0}\right) :\,\,\,\left\| f\right\| _{II,\mu
}<\infty \right\} \text{,}
\end{equation*}
where $\left\| f\right\| _{II,\mu }$ is defined by 
\begin{equation*}
\left\| f\right\| _{II,\mu }=\sup \left\{ \int \left| f\left( x\right)
j\left( d\left( x\right) \right) k\left( r\left( x\right) \right) \right|
d\nu _{0}\left( x\right) ,\,\int \left| j\right| ^{2}d\mu =\int \left|
k\right| ^{2}d\mu =1\right\} \text{.}
\end{equation*}

If $\mu _{1}$ and $\mu _{2}$ are two equivalent quasi-invariant measures,
then $\left\| f\right\| _{II,\mu _{1}}=\left\| f\right\| _{II,\mu _{2}}$,
because $\left\| f\right\| _{II,\mu }=\left\| II_{\mu }\left( \left|
f\right| \right) \right\| $ for each quasi-invariant measure $\mu $, where $%
II_{\mu }$ is the one dimensional trivial representation on $\mu $ \ 

Define $\left\| f\right\| _{II}$ to be 
\begin{equation*}
\sup \left\{ \left\| f\right\| _{II,\mu }\text{:\thinspace }\mu \text{
quasi-invariant Radon measure on }G^{\left( 0\right) }\right\}
\end{equation*}
The supremum can be taken over the classes of quasi-invariant measure.

If $\left\| {}\right\| $ is the full $C^{\ast }$-norm on $C_{c}\left(
G\right) $, then 
\begin{equation*}
\left\| f\right\| \leq \left\| f\right\| _{II}\text{ for all }f\text{.}
\end{equation*}
(see \cite{rawm})
\end{notation}

\bigskip

\begin{lemma}
\label{ll}Let $G$ be a locally compact second countable groupoid with proper
orbit space. Let $\left\{ \nu ^{u},u\in G^{\left( 0\right) }\right\} $ be a
Haar system on $G$, $\left( \left\{ \beta _{v}^{u}\right\} ,\left\{ \mu ^{%
\dot{u}}\right\} \right) $ its decomposition over the principal groupoid
associated to $G$ and $\delta $ the associated \ $1$-cocycle. If $f$ is a
universally measurable function on $G$, then 
\begin{equation*}
\left\| f\right\| _{II}\leq \sup\limits_{\dot{w}}\left( \int \int \left(
\int \left| f\left( x\right) \right| \delta \left( x\right) ^{-\frac{1}{2}%
}d\beta _{v}^{u}\left( x\right) \right) ^{2}d\mu ^{\dot{w}}\left( v\right)
d\mu ^{\dot{w}}\left( u\right) \right) ^{\frac{1}{2}}.
\end{equation*}
\end{lemma}

\begin{proof}
Each Radon quasi-invariant measure is equivalent with a Radon measure of the
form $\int \mu ^{\dot{u}}d\tilde{\mu}\left( \dot{u}\right) $, where $\tilde{%
\mu}$ is a probability measure on the orbit space $G/G^{\left( 0\right) }.$
Therefore for the computation of $\left\| .\right\| _{II}$ it is enough to
consider only the quasi-invariant measures of the form $\mu =\int \mu ^{\dot{%
u}}d\tilde{\mu}\left( \dot{u}\right) $, where $\tilde{\mu}$ is a probability
measure on $G^{\left( 0\right) }/G$. It is easy to see that the modular
function of $\int \mu ^{\dot{u}}d\tilde{\mu}\left( \dot{u}\right) $ is $%
\Delta =\delta $.

Let $j,k\in L^{2}\left( G^{\left( 0\right) },\mu \right) $ with $\int
|j|^{2}d\mu =\int |k|^{2}d\mu =1$. We have 
\begin{eqnarray*}
&&\int \int \int \int \left| f\left( x\right) \right| \delta \left( x\right)
^{-\frac{1}{2}}d\beta _{v}^{u}\left( x\right) \left| j\left( v\right)
\right| \left| k\left( u\right) \right| d\mu ^{\dot{w}}\left( v\right) d\mu
^{\dot{w}}\left( u\right) d\tilde{\mu}\left( \dot{w}\right) \\
&\leq &\int \left( \int \int \left( \int \left| f\left( x\right) \right|
\delta \left( x\right) ^{-\frac{1}{2}}d\beta _{v}^{u}\left( x\right) \right)
^{2}d\mu ^{\dot{w}}\left( v\right) d\mu ^{\dot{w}}\left( u\right) \right) ^{%
\frac{1}{2}}\cdot \\
&&\,\,\ \ \;\;\;\;\;\;\cdot \left( \int \int \left| j\left( v\right) \right|
^{2}\left| k\left( u\right) \right| ^{2}d\mu ^{\dot{w}}\left( v\right) d\mu
^{\dot{w}}\left( u\right) \right) ^{\frac{1}{2}}d\tilde{\mu}\left( \dot{w}%
\right) \\
&=&\sup\limits_{\dot{w}}\left( \int \int \left( \int \left| f\left( x\right)
\right| \delta \left( x\right) ^{-\frac{1}{2}}d\beta _{v}^{u}\left( x\right)
\right) ^{2}d\mu ^{\dot{w}}\left( v\right) d\mu ^{\dot{w}}\left( u\right)
\right) ^{\frac{1}{2}}\cdot \\
&&\,\;\;\;\ \ \ \cdot \int \left( \int \left| j\left( v\right) \right|
^{2}d\mu ^{\dot{w}}\left( v\right) \right) ^{\frac{1}{2}}\left( \left|
k\left( u\right) \right| ^{2}d\mu ^{\dot{w}}\left( u\right) \right) ^{\frac{1%
}{2}}d\tilde{\mu}\left( \dot{w}\right) \\
&\leq &\sup\limits_{\dot{w}}\left( \int \int \left( \int \left| f\left(
x\right) \right| \delta \left( x\right) ^{-\frac{1}{2}}d\beta _{v}^{u}\left(
x\right) \right) ^{2}d\mu ^{\dot{w}}\left( v\right) d\mu ^{\dot{w}}\left(
u\right) \right) ^{\frac{1}{2}}
\end{eqnarray*}
Consequently, 
\begin{equation*}
\left\| f\right\| _{II}\leq \sup\limits_{\dot{w}}\left( \int \int \left(
\int \left| f\left( x\right) \right| \delta \left( x\right) ^{-\frac{1}{2}%
}d\beta _{v}^{u}\left( x\right) \right) ^{2}d\mu ^{\dot{w}}\left( v\right)
d\mu ^{\dot{w}}\left( u\right) \right) ^{\frac{1}{2}}
\end{equation*}
\end{proof}

\bigskip

If $G$ is locally transitive, each orbit $\left[ u\right] $ is open in $%
G^{\left( 0\right) }$. Each measure $\mu ^{\dot{u}}$ is supported on $\left[
u\right] $. Since $\left( \left[ u\right] \right) $ is a partition of $%
G^{\left( 0\right) }$ into open sets, it follows that there is a unique
Radon measure $m$ on $G^{\left( 0\right) }$ such that the restriction of $m$
at $C_{c}\left( \left[ u\right] \right) $ is $\mu ^{\dot{u}}$ for each $%
\left[ u\right] $.

\begin{corollary}
\label{c}Let $G$ be a locally compact second countable \emph{locally
transitive} groupoid endowed with a Haar system $\nu =\left\{ \nu ^{u},u\in
G^{\left( 0\right) }\right\} $. Let $f$ be a universally measurable function
such that $\left\| f\right\| _{II}<\infty $.

\begin{enumerate}
\item  If $\left( f_{n}\right) _{n}$is a uniformly bounded sequence of
universally measurable functions supported on a compact set, and if $\left(
f_{n}\right) _{n}$ converges pointwise to $f$, then $\left( f_{n}\right)
_{n} $ converges to $f$ in the norm of $C^{\ast }\left( G,\nu \right) .$

\item  If $\left( f_{n}\right) _{n}$is an increasing sequence of universally
measurable nonnegative functions on $G$ that converges pointwise to $f$,
then $\left( f_{n}\right) _{n}$ converges to $f$ in the norm of $C^{\ast
}\left( G,\nu \right) .$
\end{enumerate}
\end{corollary}

\begin{proof}
Let $\left( \left\{ \beta _{v}^{u}\right\} ,\left\{ \mu ^{\dot{u}}\right\}
\right) $ be the decomposition of the Haar system over the principal
groupoid associated to $G$ and $\delta $ the associated \ $1$-cocycle. Let $%
m $ be the unique measure such that restriction of $m$ at $C_{c}\left( \left[
u\right] \right) $ is $\mu ^{\dot{u}}$ for each $\left[ u\right] $. Let $%
\left( f_{n}\right) _{n}$be a sequence of universally measurable functions
supported on a compact set $K$. \ Let us write 
\begin{equation*}
M=\sup\limits_{u,v}\beta _{u}^{v}\left( K^{-1}\right) \text{.}
\end{equation*}
and let us assume that $\left( f_{n}\right) _{n}$ converges pointwise to $f$%
. According to Lemma \ref{ll}, we have 
\begin{equation*}
\left\| f-f_{n}\right\| _{II}\leq \sup\limits_{\dot{w}}\left( \int \int
\left( \int \left| f\left( x\right) -f_{n}\left( x\right) \right| \delta
\left( x\right) ^{-\frac{1}{2}}d\beta _{v}^{u}\left( x\right) \right)
^{2}d\mu ^{\dot{w}}\left( v\right) d\mu ^{\dot{w}}\left( u\right) \right) ^{%
\frac{1}{2}}\text{.}
\end{equation*}

Hence 
\begin{eqnarray*}
\left\| f-f_{n}\right\| _{II} &\leq &\sup\limits_{\dot{w}}M\left( \int \int
\left( \int \left| f\left( x\right) -f_{n}\left( x\right) \right| ^{2}d\beta
_{v}^{u}\left( x\right) \right) d\mu ^{\dot{w}}\left( v\right) d\mu ^{\dot{w}%
}\left( u\right) \right) ^{\frac{1}{2}} \\
&\leq &M\left( \int \int \left( \int \left| f\left( x\right) -f_{n}\left(
x\right) \right| ^{2}d\beta _{u}^{v}\left( x\right) \right) dm\left(
v\right) dm\left( u\right) \right) ^{\frac{1}{2}}\text{.}
\end{eqnarray*}
If $\left\| {}\right\| $ denotes the $C^{\ast }$-norm, then 
\begin{equation*}
\lim_{n}\left\| f-f_{n}\right\| \leq \lim_{n}\left\| f-f_{n}\right\| _{II}=0%
\text{,}
\end{equation*}
because 
\begin{equation*}
\int \int \left( \int \left| f\left( x\right) -f_{n}\left( x\right) \right|
^{2}d\beta _{u}^{v}\left( x\right) \right) dm\left( v\right) dm\left(
u\right)
\end{equation*}
converges to zero, by the Dominated Convergence Theorem.

Let $\left( f_{n}\right) _{n}$be an increasing sequence of universally
measurable nonnegative functions that converges pointwise to $f.$ Since 
\begin{eqnarray*}
\left\| f-f_{n}\right\| _{II} &\leq &\sup\limits_{\dot{w}}\left( \int \int
\left( \int \left| f\left( x\right) -f_{n}\left( x\right) \right| \delta
\left( x\right) ^{-\frac{1}{2}}d\beta _{v}^{u}\left( x\right) \right)
^{2}d\mu ^{\dot{w}}\left( v\right) d\mu ^{\dot{w}}\left( u\right) \right) ^{%
\frac{1}{2}} \\
&\leq &\left( \int \int \left( \int \left| f\left( x\right) -f_{n}\left(
x\right) \right| \delta \left( x\right) ^{-\frac{1}{2}}d\beta _{v}^{u}\left(
x\right) \right) ^{2}dm\left( v\right) dm\left( u\right) \right) ^{\frac{1}{2%
}}
\end{eqnarray*}
it follows that 
\begin{equation*}
\lim_{n}\left\| f-f_{n}\right\| _{II}=0\text{.}
\end{equation*}
\end{proof}

\begin{proposition}
\label{elem}Let $G$ be a locally compact second countable \emph{locally
transitive }groupoid endowed with a Haar system $\nu =\left\{ \nu ^{u},u\in
G^{\left( 0\right) }\right\} $.Then any function in $\mathcal{B}_{c}\left(
G\right) $, the space of compactly supported Borel bounded function on $G$,
can be viewed as an element of $C^{\ast }\left( G,\nu \right) $.
\end{proposition}

\begin{proof}
Let $\left( \left\{ \beta _{v}^{u}\right\} ,\left\{ \mu ^{\dot{u}}\right\}
\right) $ be the decomposition of the Haar system over the principal
groupoid associated to $G$, and $\delta $ the associated \ $1$-cocycle. Let $%
m$ be the unique measure such that restriction of $m$ at $C_{c}\left( \left[
u\right] \right) $ is $\mu ^{\dot{u}}$ for each $\left[ u\right] $. Let $m$
be a dominant for the family $\left\{ \mu ^{\dot{u}}\right\} $. Let $\nu
_{1} $ be the measure on $G$ define by 
\begin{equation*}
\int f\left( x\right) d\nu _{1}\left( x\right) =\left( \int \int \left( \int
f\left( x\right) d\beta _{u}^{v}\left( x\right) \right) dm\left( v\right)
dm\left( u\right) \right)
\end{equation*}
for all Borel nonnegative function $f$. If $f\in $ $\mathcal{B}_{c}\left(
G\right) $, then $f$ is the limit in $L^{2}\left( G,\nu _{1}\right) $ of a
sequence, $\left( f_{n}\right) _{n}$, in $C_{c}\left( G\right) $ that is
supported on some compact set $K$ supporting $f$. If we write 
\begin{equation*}
M=\sup\limits_{u,v}\beta _{u}^{v}\left( K^{-1}\right) \text{,}
\end{equation*}
we have 
\begin{eqnarray*}
\left\| f-f_{n}\right\| _{II} &\leq &\sup\limits_{\dot{w}}M\left( \int \int
\left( \int \left| f\left( x\right) -f_{n}\left( x\right) \right| ^{2}d\beta
_{v}^{u}\left( x\right) \right) d\mu ^{\dot{w}}\left( v\right) d\mu ^{\dot{w}%
}\left( u\right) \right) ^{\frac{1}{2}} \\
&\leq &M\left( \int \int \left( \int \left| f\left( x\right) -f_{n}\left(
x\right) \right| ^{2}d\beta _{u}^{v}\left( x\right) \right) dm\left(
v\right) dm\left( u\right) \right) ^{\frac{1}{2}}\text{.}
\end{eqnarray*}
If $\left\| {}\right\| $ denotes the $C^{\ast }$-norm, then 
\begin{equation*}
\lim_{n}\left\| f-f_{n}\right\| \leq \lim_{n}\left\| f-f_{n}\right\| _{II}=0%
\text{. }
\end{equation*}
Thus $f$ can be viewed as an element in $C^{\ast }\left( G,\nu \right) $.
\end{proof}

\begin{proposition}
If $G$ is a locally compact second countable \emph{locally transitive }%
groupoid endowed with a Haar system $\left\{ \nu ^{u},\,u\in G^{\left(
0\right) }\right\} $ with bounded decomposition, then 
\begin{equation*}
C^{\ast }\left( G,\nu \right) =M^{\ast }\left( G,\nu \right) \text{.}
\end{equation*}
\end{proposition}

\begin{proof}
It follows from the Proposition \ref{elem}.
\end{proof}

\bigskip

\begin{remark}
\label{giso}Let $G$ be locally compact locally transitive groupoid. Let $F$
be a subset of $G^{\left( 0\right) }$ containing only one element $e\left(
u\right) $ in each orbit $\left[ u\right] $. It is easy to see that $F$ is a
closed subset of $G$ and that $F$ is a discrete space. Let $\sigma
:G^{\left( 0\right) }\rightarrow G^{F}$ be a \textit{regular cross section}
of $d_{F}$. Let us endow $\bigcup\limits_{\left[ u\right] }\left[ u\right]
\times G_{e\left( u\right) }^{e\left( u\right) }\times \left[ u\right] $
with the topology induced from $G^{\left( 0\right) }\times G_{F}^{F}\times
G^{\left( 0\right) }$. The topology of $\bigcup\limits_{\left[ u\right] }%
\left[ u\right] \times G_{e\left( u\right) }^{e\left( u\right) }\times \left[
u\right] $ is locally compact because $\bigcup\limits_{\left[ u\right] }%
\left[ u\right] \times G_{e\left( u\right) }^{e\left( u\right) }\times \left[
u\right] $ is a closed subset of the locally compact space $G^{\left(
0\right) }\times G_{F}^{F}\times G^{\left( 0\right) }$. With the operations: 
\begin{eqnarray*}
\left( u,x,v\right) \left( v,y,w\right) &=&\left( u,xy,w\right) \\
\left( u,x,v\right) ^{-1} &=&\left( v,x^{-1},u\right)
\end{eqnarray*}
$\bigcup\limits_{\left[ u\right] }\left[ u\right] \times G_{e\left( u\right)
}^{e\left( u\right) }\times \left[ u\right] $ becomes a groupoid. Let us
define $\phi :G\rightarrow \bigcup\limits_{\left[ u\right] }\left[ u\right]
\times G_{e\left( u\right) }^{e\left( u\right) }\times \left[ u\right] $ by 
\begin{equation*}
\phi \left( x\right) =\left( r\left( x\right) ,\sigma \left( r\left(
x\right) \right) x\sigma \left( d\left( x\right) \right) ^{-1},d\left(
x\right) \right)
\end{equation*}
and note that $\phi $ is a Borel isomorphism which carries compact sets to
relatively compact sets .
\end{remark}

\begin{lemma}
Let $G$ be locally compact second countable \emph{locally transitive}
groupoid. Let $F$ be a subset of $G^{\left( 0\right) }$ containing only one
element $e\left( u\right) $ in each orbit $\left[ u\right] $. Let $\sigma
:G^{\left( 0\right) }\rightarrow G^{F}$ be a \textit{regular cross section}
of $d_{F}$. Then any compactly supported Borel bounded function on $G$ is
pointwise limit of a uniformly bounded sequence $\left( f_{n}\right) _{n}$
of Borel functions supported on a compact set supporting $f$, having the
property that each $f_{n}$ is a linear combination of functions of the form 
\begin{equation*}
x\rightarrow g_{1}\left( r\left( x\right) \right) g\left( \sigma \left(
r\left( x\right) \right) x\sigma \left( d\left( x\right) \right)
^{-1}\right) g_{2}\left( d\left( x\right) \right)
\end{equation*}
where $g_{1},g_{2}$ are compactly supported Borel bounded function on $%
G^{\left( 0\right) }$ and $g$ is a compactly supported Borel bounded
function on $G_{F}^{F}$.
\end{lemma}

\begin{proof}
Let us endow $\bigcup\limits_{\left[ u\right] }\left[ u\right] \times
G_{e\left( u\right) }^{e\left( u\right) }\times \left[ u\right] $ with the
topology induced from $G^{\left( 0\right) }\times G_{F}^{F}\times G^{\left(
0\right) }$ as in Remark.\ref{giso}. The topology of $\bigcup\limits_{\left[
u\right] }\left[ u\right] \times G_{e\left( u\right) }^{e\left( u\right)
}\times \left[ u\right] $ is locally compact. Any compactly supported Borel
bounded function on $G^{\left( 0\right) }\times G_{F}^{F}\times G^{\left(
0\right) }$ is pointwise limit of uniformly bounded sequences $\left(
f_{n}\right) _{n}$ of Borel functions supported on a compact set, such that
each function $f_{n}$ is a linear combination of functions of the form 
\begin{equation*}
\left( u,x,v\right) \rightarrow g_{1}\left( u\right) g\left( x\right)
g_{2}\left( v\right)
\end{equation*}
where $g_{1},g_{2}$ are compactly supported Borel bounded function on $%
G^{\left( 0\right) }$ and $g$ is a compactly supported Borel bounded
function on $G_{F}^{F}$ . Consequently, any compactly supported Borel
bounded function on $\bigcup\limits_{\left[ u\right] }\left[ u\right] \times
G_{e\left( u\right) }^{e\left( u\right) }\times \left[ u\right] $ has the
same property. Since $\phi :G\rightarrow \bigcup\limits_{\left[ u\right] }%
\left[ u\right] \times G_{e\left( u\right) }^{e\left( u\right) }\times \left[
u\right] $ defined by 
\begin{equation*}
\phi \left( x\right) =\left( r\left( x\right) ,\sigma \left( r\left(
x\right) \right) x\sigma \left( d\left( x\right) \right) ^{-1},d\left(
x\right) \right)
\end{equation*}
is a Borel isomorphism which carries compact sets to relatively compact
sets, it follows that any compactly supported Borel bounded function on $G$
can be represented as a pointwise limit of a uniformly bounded sequence $%
\left( f_{n}\right) _{n}$ of Borel functions supported on a compact set
supporting $f$ , having the property that each $f_{n}$ is a linear
combination of functions of the form 
\begin{equation*}
x\rightarrow g_{1}\left( r\left( x\right) \right) g\left( \sigma \left(
r\left( x\right) \right) x\sigma \left( d\left( x\right) \right)
^{-1}\right) g_{2}\left( d\left( x\right) \right)
\end{equation*}
\end{proof}

\begin{corollary}
Let $G$ be locally compact second countable \emph{locally transitive}
groupoid. Let $F$ be a subset of $G^{\left( 0\right) }$ containing only one
element $e\left( u\right) $ in each orbit $\left[ u\right] $. Let $\sigma
:G^{\left( 0\right) }\rightarrow G^{F}$ be a \textit{regular cross section}
of $d_{F}$. Then the linear span of the functions of the form 
\begin{equation*}
x\rightarrow g_{1}\left( r\left( x\right) \right) g\left( \sigma \left(
r\left( x\right) \right) x\sigma \left( d\left( x\right) \right)
^{-1}\right) g_{2}\left( d\left( x\right) \right)
\end{equation*}
where $g_{1},g_{2}\in \mathcal{B}_{c}\left( G^{\left( 0\right) }\right) $
and $g\in \mathcal{B}_{c}\left( G_{F}^{F}\right) ,$ is dense in the full $%
C^{\ast }$-algebra of $G$.
\end{corollary}

\begin{proof}
Let $f$ be a function on $G$, defined by 
\begin{equation*}
f\left( x\right) =g_{1}\left( r\left( x\right) \right) g\left( \sigma \left(
r\left( x\right) \right) x\sigma \left( d\left( x\right) \right)
^{-1}\right) g_{2}\left( d\left( x\right) \right)
\end{equation*}
where $g_{1},g_{2}\in \mathcal{B}_{c}\left( G^{\left( 0\right) }\right) $
and $g\in \mathcal{B}_{c}\left( G_{F}^{F}\right) $. Then $f\in \mathcal{B}%
_{c}\left( G\right) $, therefore it can be viewed as an element of the $%
C^{\ast }\left( G,\nu \right) $ as we note in Proposition \ref{elem}. Each $%
f\in \mathcal{B}_{c}\left( G\right) $ (in particular in $C_{c}\left(
G\right) $) is the limit (pointwise and consequently in the $C^{\ast }$-norm
according to Corollary \ref{c}) of a uniformly bounded sequence $\left(
f_{n}\right) _{n}$ of Borel functions supported on a compact set supporting $%
f$, having the property that each $f_{n}$ is a linear combination of
functions of the required form.
\end{proof}

\bigskip \bigskip

\begin{proposition}
Let $G$ is a locally compact second countable \emph{locally transitive}
groupoid endowed with a Haar system $\left\{ \nu ^{u},\,u\in G^{\left(
0\right) }\right\} $. Let $F$ be a subset of $G^{\left( 0\right) }$
containing only one element $e\left( u\right) $ in each orbit $\left[ u%
\right] $. Let $\sigma :G^{\left( 0\right) }\rightarrow G^{F}$ be a \textit{%
regular cross section} of $d_{F}$. Then 
\begin{equation*}
C^{\ast }\left( G,\nu \right) =M^{\ast }\left( G,\nu \right) =M_{\sigma
}^{\ast }\left( G,\nu \right) \text{.}
\end{equation*}
\end{proposition}

\begin{proof}
We have proved that \bigskip $C^{\ast }\left( G,\nu \right) =M^{\ast }\left(
G,\nu \right) $. From the preceding corollary, it follows that the linear
span of the functions of the form 
\begin{equation*}
x\rightarrow g_{1}\left( r\left( x\right) \right) g\left( \sigma \left(
r\left( x\right) \right) x\sigma \left( d\left( x\right) \right)
^{-1}\right) g_{2}\left( d\left( x\right) \right)
\end{equation*}
where $g_{1},g_{2}\in \mathcal{B}_{c}\left( G^{\left( 0\right) }\right) $
and $g\in \mathcal{B}_{c}\left( G_{F}^{F}\right) $ is dense in\ $C^{\ast
}\left( G,\nu \right) $. But this space is contained in $\mathcal{B}_{\sigma
}\left( G\right) $. Therefore $C^{\ast }\left( G,\nu \right) =M^{\ast
}\left( G,\nu \right) =M_{\sigma }^{\ast }\left( G,\nu \right) $.
\end{proof}

\bigskip \bigskip

\section{The case of principal proper groupoids case}

\begin{notation}
\label{topo}Let $G$ be a locally compact second countable groupoid with
proper orbit space. Let $F$ be a Borel subset of $G^{\left( 0\right) }$
containing only one element $e\left( u\right) $ in each orbit $\left[ u%
\right] $. Let $\sigma :G^{\left( 0\right) }\rightarrow G^{F}$ be a\ regular
cross section for $d_{F}:G^{F}\rightarrow G^{\left( 0\right) },d_{F}\left(
x\right) =d\left( x\right) $ with $\sigma \left( e\left( v\right) \right)
=e\left( v\right) $ for all $v\in G^{\left( 0\right) }$. Let $q:G\rightarrow
G_{F}^{F}$ \ be defined by 
\begin{equation*}
q\left( x\right) =\sigma \left( r\left( x\right) \right) x\sigma \left(
d\left( x\right) \right) ^{-1}
\end{equation*}
\bigskip We shall endow $G_{F}^{F}$ \ with the quotient topology induced by $%
q.$ We shall denote by $\mathbf{C}_{\sigma }\left( G\right) $ the linear
span of the functions of the form 
\begin{equation*}
x\rightarrow g_{1}\left( r\left( x\right) \right) g\left( \sigma \left(
r\left( x\right) \right) x\sigma \left( d\left( x\right) \right)
^{-1}\right) g_{2}\left( d\left( x\right) \right)
\end{equation*}
where $g_{1},g_{2}\in C_{c}\left( G^{\left( 0\right) }\right) $ and $g\in
C_{c}\left( G_{F}^{F}\right) $.
\end{notation}

\begin{proposition}
\label{dens}With the Notation \ref{topo}, if \ the space of continuous
functions with compact support on $G_{F}^{F}$ (with the respect to the
quotient topology induced by $q$) separates the points of $G_{F}^{F}$, then $%
\mathbf{C}_{\sigma }\left( G\right) $ is dense in $C_{c}\left( G\right) $
(for the inductive limit topology). In particular, if the quotient topology
induced by $q$ on $G_{F}^{F}$ is a locally compact (Hausdorff) topology,
then $\mathbf{C}_{\sigma }\left( G\right) $ is dense in $C_{c}\left(
G\right) $.
\end{proposition}

\begin{proof}
If \ the space of continuous functions with compact support on $G_{F}^{F}$
(with the respect to the quotient topology induced by $q$) separates the
points of $G_{F}^{F}$, then $\mathbf{C}_{\sigma }\left( G\right) $ separates
the points of $G$. By Stone-Weierstrass Theorem, it follows that $\mathbf{C}%
_{\sigma }\left( G\right) $ is dense in $C_{c}\left( G\right) $ (for the
inductive limit topology)
\end{proof}

\begin{proposition}
\label{princ}Let $G$ be a locally compact principal groupoid. If $G$ is
proper, then the quotient topology induced by $q$ on $G_{F}^{F}$ is a
locally compact (Hausdorff) topology. \ Consequently, $\mathbf{C}_{\sigma
}\left( G\right) $ is dense in $C_{c}\left( G\right) $ for the inductive
limit topology (we use the Notation \ref{topo}).
\end{proposition}

\begin{proof}
Let $\pi :G\rightarrow G^{\left( 0\right) }/G$ be the canonical projection.
Let us note that for a principal groupoid the condition 
\begin{equation*}
q\left( x\right) =q\left( y\right)
\end{equation*}
is equivalent with 
\begin{equation*}
\pi \left( r\left( x\right) \right) =\pi \left( r\left( y\right) \right) 
\text{.}
\end{equation*}

First we shall prove that the topology on $G_{F}^{F}$ is Hausdorff. Let $%
\left( x_{i}\right) _{i}$ and $\left( y_{i}\right) _{i}$ be two nets with $%
q\left( x_{i}\right) =q\left( y_{i}\right) $ for every $i$. Let us suppose
that $\left( x_{i}\right) _{i}$ converges to $x$ and $\left( y_{i}\right)
_{i}$ converges to $y$. Then 
\begin{equation*}
\lim \pi \left( r\left( x_{i}\right) \right) =\lim \pi \left( r\left(
y_{i}\right) \right) =\pi \left( r\left( x\right) \right) =\pi \left(
r\left( y\right) \right)
\end{equation*}

Hence $q\left( x\right) =q\left( y\right) $, and therefore the topology on $%
G_{F}^{F}$ is Hausdorff. We shall prove that $q$ is open. If $\left(
z_{i}\right) _{i}$ is a net converging to $q\left( x\right) $ in $G_{F}^{F}$%
, then $\pi \circ r\left( z_{i}\right) $ converges to $\pi \circ r\left(
x\right) $ . Since 
\begin{equation*}
\pi \circ r:G\rightarrow G^{\left( 0\right) }/G
\end{equation*}
is an open map, there is a net $\left( x_{i}\right) _{i}$ converging to $x$,
such that $\pi \circ r\left( x_{i}\right) =\pi \circ r\left( z_{i}\right) $,
and consequently $q\left( x_{i}\right) =q\left( z_{i}\right) =z_{i}$. \
Hence $q$ is an open map and the quotient topology induced by $q$ on $%
G_{F}^{F}$ is locally compact.
\end{proof}

\begin{theorem}
\label{con}Let $G$ be a locally compact second countable groupoid with
proper orbit space. Let $F$ be a Borel subset of $G^{\left( 0\right) }$
meeting each orbit exactly once. Let $\sigma :G^{\left( 0\right)
}\rightarrow G^{F}$ be a\textit{\ }regular cross section for $%
d:G^{F}\rightarrow G$. Let us assume that the quotient topology induced by $%
q $ on $G_{F}^{F}$ is a locally compact (Hausdorff) topology. Let $\left\{
\nu ^{u},u\in G^{\left( 0\right) }\right\} $ be a Haar system on $G$. Then 
\begin{equation*}
C^{\ast }\left( G,\nu \right) \subset M_{\sigma }^{\ast }\left( G,\nu
\right) \subset M^{\ast }\left( G,\nu \right) \text{.}
\end{equation*}
\end{theorem}

\begin{proof}
From Proposition \ref{dens} $\mathbf{C}_{\sigma }\left( G\right) $ is dense
in $C_{c}\left( G\right) $ for the inductive limit topology and hence is
dense in $C^{\ast }\left( G,\nu \right) $ . Since $\mathbf{C}_{\sigma
}\left( G\right) \subset \mathcal{B}_{\sigma }\left( G\right) $, it follows
that $C^{\ast }\left( G,\nu \right) \subset M_{\sigma }^{\ast }\left( G,\nu
\right) $.
\end{proof}

\begin{corollary}
Let $G$ be a locally compact second countable\emph{\ principal proper}
groupoid. Let $F$ be a Borel subset of $G^{\left( 0\right) }$ meeting each
orbit exactly once. Let $\sigma :G^{\left( 0\right) }\rightarrow G^{F}$ be a%
\textit{\ }regular cross section for $d:G^{F}\rightarrow G$. Let $\left\{
\nu ^{u},u\in G^{\left( 0\right) }\right\} $ be a Haar system on $G$. Then 
\begin{equation*}
C^{\ast }\left( G,\nu \right) \subset M_{\sigma }^{\ast }\left( G,\nu
\right) \subset M^{\ast }\left( G,\nu \right) \text{.}
\end{equation*}
\end{corollary}

\begin{proof}
Applying Proposition \ref{princ}, we obtain that the quotient topology
induced by $q$ on $G_{F}^{F}$ is a locally compact (Hausdorff) topology.
Therefore $G$ satisfies the hypothesis of Theorem \ref{con}.
\end{proof}

\bigskip

\begin{definition}
\label{bcont}Let $\left\{ \mu ^{\dot{u}}\right\} _{\dot{u}}$ be a system of
measures on $G^{\left( 0\right) }$ satisfying:

\begin{enumerate}
\item  $supp\left( \mu ^{\dot{u}}\right) =\left[ u\right] $ for all $\dot{u}$%
.

\item  For all compactly supported continuous functions $f$ on $G^{\left(
0\right) }$ the function 
\begin{equation*}
u\rightarrow \int f\left( v\right) \mu ^{\pi \left( u\right) }\left( v\right)
\end{equation*}
is continuous
\end{enumerate}

We shall say that the Hilbert bundle determined by the system of measures $%
\left\{ \mu ^{\dot{u}}\right\} _{\dot{u}}$ has a continuous basis if there
is sequence $f_{1},f_{2},...$ of real valued \emph{continuous} functions on $%
G^{\left( 0\right) }$ such that $\dim \left( L^{2}\left( \mu ^{\dot{u}%
}\right) \right) =\infty $ if and only if $\left\| f_{n}\right\| _{2}$ for $%
n=1,2,...$ and then $\left\{ f_{1},f_{2},...\right\} $ gives an ortonormate
basis of $L^{2}\left( \mu ^{\dot{u}}\right) $, while $\dim \left(
L^{2}\left( \mu ^{\dot{u}}\right) \right) =k<\infty $ if and only if $%
\left\| f_{n}\right\| _{2}=1$ for $n\leq k$, and $\left\| f_{n}\right\|
_{2}=0$ for $n>k$ \ and then $\left\{ f_{1},f_{2},...f_{k}\right\} $ gives
an ortonormate basis of $L^{2}\left( \mu ^{\dot{u}}\right) $.
\end{definition}

\begin{remark}
\label{r1}Let $\left\{ \mu _{1}^{\dot{u}}\right\} _{\dot{u}}$ and $\left\{
\mu _{2}^{\dot{u}}\right\} _{\dot{u}}$ be two systems of measures on $%
G^{\left( 0\right) }$ satisfying:\bigskip

\begin{enumerate}
\item  $supp\left( \mu _{i}^{\dot{u}}\right) =\left[ u\right] $ for all $%
\dot{u}$, $i=1,2$.

\item  For all compactly supported continuous functions $f$ on $G^{\left(
0\right) }$ the function 
\begin{equation*}
u\rightarrow \int f\left( v\right) \mu _{i}^{\pi \left( u\right) }\left(
v\right)
\end{equation*}
is continuous
\end{enumerate}

Let us assume that the Hilbert bundles determined by the systems of measures 
$\left\{ \mu _{i}^{\dot{u}}\right\} _{\dot{u}}$ have continuous bases. Let \ 
$f_{1}.f_{2},...$ be a continuous basis for Hilbert bundle determined by $%
\left\{ \mu _{1}^{\dot{u}}\right\} _{\dot{u}}$ and let $g_{1},g_{2},...$ be
a continuous basis for Hilbert bundle determined by $\left\{ \mu _{2}^{\dot{u%
}}\right\} _{\dot{u}}$. Let us define a unitary operator $U_{\dot{u}%
}:L^{2}\left( \mu _{1}^{\dot{u}}\right) \rightarrow L^{2}\left( \mu _{2}^{%
\dot{u}}\right) $ by 
\begin{equation*}
U_{\dot{u}}\left( f_{n}\right) =g_{n}\text{ for all }n
\end{equation*}
Then the family $\left\{ U_{\dot{u}}\right\} _{\dot{u}}$ has the following
properties:

\begin{enumerate}
\item  For all Borel bounded function $f$ on $G^{\left( 0\right) }$, 
\begin{equation*}
u\rightarrow U_{\pi \left( u\right) }\left( f\right)
\end{equation*}
is a bounded Borel function with compact support.

\item  For all Borel bounded function $f$ on $G^{\left( 0\right) }$, 
\begin{equation*}
U_{\pi \left( u\right) }\left( \overline{f}\right) =\overline{U_{\pi \left(
u\right) }\left( f\right) }
\end{equation*}

\item  For all compactly supported continuous functions $f$ on $G^{\left(
0\right) }$ there is a sequence $\left( h_{n}\right) _{n}$ of compactly
supported continuous functions on $G^{\left( 0\right) }$ such that 
\begin{equation*}
sup_{\dot{u}}\int \left| U_{\dot{u}}\left( f\right) -h_{n}\right| ^{2}d\mu
_{2}^{\dot{u}}\rightarrow 0\,\left( n\rightarrow \infty \right)
\end{equation*}
\end{enumerate}

Indeed, we can define 
\begin{equation*}
h_{n}\left( v\right) =\sum_{k=1}^{n}g_{k}\left( v\right) \int f\left(
u\right) f_{k}\left( u\right) d\mu _{1}^{\pi \left( v\right) }\left(
u\right) \text{.}
\end{equation*}
\end{remark}

\bigskip

\begin{remark}
\label{r2}Let $G$ be a locally compact second countable groupoid with proper
orbit space. Let $F$ be a Borel subset of $G^{\left( 0\right) }$ containing
only one element $e\left( u\right) $ in each orbit $\left[ u\right] $. Let
us assume that $F\cap \left[ K\right] $ has a compact closure for each
compact subset $K$ of $G^{\left( 0\right) }$, and let $\sigma :G^{\left(
0\right) }\rightarrow G^{F}$ be a\ regular cross section for $%
d_{F}:G^{F}\rightarrow G^{\left( 0\right) }$. Let us endow $\ G_{F}^{F}$ \
with the quotient topology induced by $q:G\rightarrow G_{F}^{F}$ \ 
\begin{equation*}
q\left( x\right) =\sigma \left( r\left( x\right) \right) x\sigma \left(
d\left( x\right) \right) ^{-1},\,x\in G
\end{equation*}
\bigskip If $g\in C_{c}\left( G_{F}^{F}\right) $ and $g_{1},\,g_{2}$ are two
functions on $G^{\left( 0\right) }$ with the property that there is two
sequences $\left( h_{n}^{1}\right) _{n}$ and $\left( h_{n}^{2}\right) _{n}$
of compactly supported continuous functions on $G^{\left( 0\right) }$ such
that 
\begin{equation*}
sup_{\dot{u}}\int \left| g_{i}-h_{n}^{i}\right| ^{2}d\mu _{2}^{\dot{u}%
}\rightarrow 0\,\left( n\rightarrow \infty \right)
\end{equation*}
for $i=1,2$, then

\begin{equation*}
x\overset{f}{\rightarrow }g_{1}\left( r\left( x\right) \right) g\left(
\sigma \left( r\left( x\right) \right) x\sigma \left( d\left( x\right)
\right) ^{-1}\right) g_{2}\left( d\left( x\right) \right)
\end{equation*}
can be viewed as an element of $C^{\ast }\left( G,\nu \right) $. Indeed, it
is easy to see that 
\begin{equation*}
\left\| f-\left( h_{n}^{1}\circ r\right) \left( g\circ q\right) \left(
h_{n}^{2}\circ d\right) \right\| _{II}\rightarrow 0\,\,\left( n\rightarrow
\infty \right) .
\end{equation*}
\end{remark}

\begin{proposition}
Let $G$ be a locally compact second countable\emph{\ principal proper}
groupoid. Let $\nu _{i}=\left\{ \nu _{i}^{u},u\in G^{\left( 0\right)
}\right\} ,$ $i=1,2$ be two Haar system on $G$ and $\left( \left\{ \beta
_{v}^{u}\right\} ,\left\{ \mu _{i}^{\dot{u}}\right\} \right) $ the
corresponding decompositions over the principal groupoid. If the Hilbert
bundles determined by the systems of measures $\left\{ \mu _{i}^{\dot{u}%
}\right\} _{\dot{u}}$ have continuous bases, then the $C^{\ast }$-algebras $%
C^{\ast }\left( G,\nu _{1}\right) $ and $C^{\ast }\left( G,\nu _{2}\right) $
are $\ast $-isomorphic.
\end{proposition}

\begin{proof}
We use the Notation \ref{topo}. From Proposition \ref{dens} $\mathbf{C}%
_{\sigma }\left( G\right) $ is dense in $C_{c}\left( G\right) $ for the
inductive limit topology and hence is dense in $C^{\ast }\left( G,\nu
_{1}\right) $. We shall define a $\ast $-homomorphism $\Phi $ from $\mathbf{C%
}_{\sigma }\left( G\right) $ to $C^{\ast }\left( G,\nu _{2}\right) $. It
suffices to define $\Phi $ on the set of function on $G$ of the form 
\begin{equation*}
x\rightarrow g_{1}\left( r\left( x\right) \right) g\left( q\left( x\right)
\right) g_{2}\left( d\left( x\right) \right)
\end{equation*}
where $g_{1},g_{2}\in C_{c}\left( G^{\left( 0\right) }\right) $ and $g\in
C_{c}\left( G_{F}^{F}\right) $. Let $\left\{ U_{\dot{u}}\right\} _{\dot{u}}$
be the family of unitary operators with the properties stated in Remark \ref
{r1} associated to the systems of measures $\left\{ \mu _{i}^{\dot{u}%
}\right\} _{\dot{u}}$, $i=1,2$.

Let us define $\Phi $ by 
\begin{equation*}
\Phi \left( f\right) =\left( x\rightarrow U_{\pi \left( r\left( x\right)
\right) }\left( g_{1}\right) \left( r\left( x\right) \right) g\left( q\left(
x\right) \right) U_{\pi \left( d\left( x\right) \right) }\left( g_{2}\right)
\left( d\left( x\right) \right) \right)
\end{equation*}
where $f$ is defined by 
\begin{equation*}
f\left( x\right) =g_{1}\left( r\left( x\right) \right) g\left( q\left(
x\right) \right) g_{2}\left( d\left( x\right) \right)
\end{equation*}
with $g_{1},g_{2}\in C_{c}\left( G^{\left( 0\right) }\right) $ and $g\in
C_{c}\left( G_{F}^{F}\right) $.

We noted in Remark \ref{r2} that the functions of the form $\Phi \left(
f\right) $ can be viewed as elements of $C^{\ast }\left( G,\nu _{2}\right) $%
. With the same argument as in the proof of Theorem \ref{iso}, it follows
that $\Phi $ can be extended to $\ast $-isomorphism between $C^{\ast }\left(
G,\nu _{1}\right) $ and $C^{\ast }\left( G,\nu _{2}\right) $.
\end{proof}

\bigskip

University Constantin Br\^{a}ncu\c{s}i of T\^{a}rgu-Jiu

Bulevardul Republicii, Nr. 1,

210152 T\^{a}rgu-Jiu , Gorj

Romania

e-mail: ada@utgjiu.ro

\end{document}